\documentclass{amsart}
\usepackage{latexsym}
\usepackage{amsmath}
\usepackage{amsfonts}
\usepackage{amssymb}
\usepackage{graphicx}
\newtheorem{thm}{Theorem}[section]

\newtheorem{lem}[thm]{Lemma}
\newtheorem{prop}[thm]{Proposition}
\theoremstyle{definition}

\theoremstyle{remark}

\numberwithin{equation}{section}

\begin{document}

\title{Nonlinear Instability in Gravitational Euler-Poisson system
for $\gamma=\frac{6}{5}$}
\author{Juhi Jang}
\address{Department of Mathematics,
 Brown University, Providence,
RI 02912, USA}
\email{juhijang@math.brown.edu}


\begin{abstract}
The dynamics of gaseous stars can be described by the Euler-Poisson
system. Inspired by Rein's stability result for $\gamma>\frac{4}{3}$,
we prove the nonlinear instability of steady states for the
adiabatic exponent $\gamma=\frac{6}{5}$ in spherically symmetric and
isentropic motion.
\end{abstract}

\maketitle
\section{Introduction}
The motion of gaseous stars can be described by the Euler-Poisson
equations:

\begin{equation}
\frac{\partial\rho}{\partial t} + \nabla\cdot(\rho u)=0
\end{equation}
\begin{equation}
\frac{\partial(\rho u)}{\partial t} + \nabla\cdot(\rho u \otimes
u) + \nabla p = -\rho\nabla\Phi
\end{equation}
\begin{equation}
\triangle \Phi = 4\pi \rho
\
\end{equation}
where $t\geq 0$, $x\in\mathbb{R}^3$, $\rho$ is the density,
$u\in\mathbb{R}^3$ the velocity, $p$ the pressure of the gas, and
$\Phi$ the potential function of the self-gravitational force. We
consider the isentropic motion i.e.
$p=A\rho^\gamma$($1<\gamma<2$), where $A$ is an entropy constant
and $\gamma$ is an adiabatic exponent. In our case $A$ will be
normalized as $\frac{2\pi}{9}$ and $\gamma$ will be chosen as
$\frac{6}{5}$.\

For the spherically symmetric motion the above equations where
$\rho>0$ can be written as following:

\begin{equation}
\frac{\partial\rho}{\partial t} + u \frac{\partial\rho}{\partial
r} + \rho\frac{\partial u}{\partial r} + \frac{2}{r} \rho u = 0,
\end{equation}
\begin{equation}
\frac{\partial u}{\partial t} + u \frac{\partial u}{\partial r} +
\frac{1}{\rho}\frac{\partial p}{\partial r} +
\frac{4\pi}{r^2}\int_{0}^{r}\rho s^2 ds = 0.
\end{equation}\

First consider the stationary solutions $(\rho_0(r), 0)$ of (1.4)
and (1.5). They satisfy
\begin{equation}
\frac{dp}{dr} + \frac{4\pi\rho}{r^2}\int_{0}^{r}\rho s^2 ds = 0.
\end{equation}

 This ordinary differential equation has been well studied.
 One interesting question relevant to our context
 can be given like this: Given the mass $M>0$, how many solutions
 are there for (1.6) with $\int_{\mathbb{R}^3}\rho dx = M(\rho)= M$?
 For our purpose we summarize the answer according to the range of $\gamma$.
 See \cite{ch} or \cite{lin} for details. If $\frac{6}{5}<\gamma<2$ and any $M>0$,
 there exists at least one compactly supported stationary solution
 $\rho$ such that $M(\rho)=M$. For $\frac{4}{3}<\gamma<2$,
 every stationary solution is compactly supported and unique.
If $\gamma=\frac{6}{5}$ and any $M>0$, there is a unique
ground-state solution (not compactly supported) $\rho$. The
solution can be expressed in terms of Lane-Emden function and
moreover it can be written explicitly as
$\rho_0(r)=\frac{1}{(1+r^2)^{\frac{5}{2}}}$ up to scaling in r
and constant multiplication. On the other hand, if
$1<\gamma<\frac{6}{5}$, there is no stationary solution with
finite total mass.

The stability question has been a great interest and it has been
conjectured by astrophysicists that stationary solutions for
$\gamma<\frac{4}{3}$ are unstable; indeed one can easily check
that when $\gamma\in(1, \frac{4}{3})$ steady states are not
minimizers of the energy functional
\begin{equation}
E=E(\rho,u)=\int(\frac{1}{2}\rho |u|^2 + \frac{1}{\gamma-1}p)dx
-\frac{1}{2}\int\int\frac{\rho(x)\rho(y)}{|x-y|}dxdy
\end{equation}
by constructing a scaling invariant family of steady states, and
this indicates the possibility of certain kind of instability. So
far only partial results are known in this direction. The linear
stability of the above stationary solutions was studied in \cite{lin}
in Lagrangian formulation. It was shown that any stationary
solution is stable when $\gamma\in(\frac{4}{3},2)$ and unstable
when $\gamma\in(1, \frac{4}{3})$. In accordance with the linear
stability, a nonlinear stability for $\gamma>\frac{4}{3}$ was
established recently in \cite{r} by Rein using the variational
approach based on the fact the steady states are minimizers of the
energy functional $E$ defined in (1.7). See \cite{r2} for a
great overview of mathematical results on the nonlinear stability
problems.
For $\gamma=\frac{4}{3}$, the energy of a
steady state is zero and any small perturbation can make the
energy positive and cause part of the system go off to infinity,
which implies an instability of such a state. This kind of
instability was investigated in \cite{dlyy}. However, the
method they employed is not applicable when $\gamma<\frac{4}{3}$.
The stability question for $\frac{6}{5}\leq\gamma<\frac{4}{3}$
under the physical consideration waits for a satisfactory answer.

Our main result in this paper concerns a fully nonlinear,
dynamical instability of the steady profile for
$\gamma=\frac{6}{5}$ :
\begin{equation}
\rho_0(r)=\frac{1}{(1+r^2)^{\frac{5}{2}}}\text{ }\text{ and
}\text{ }u_0(r)=0
\end{equation}
We shall show that this steady profile is unstable under an
appropriate energy-like measurement $\mathcal{E}_l$ which will be
precisely defined later in Section 5. Here $l\leq 0$ represents
the strength of weights. As $l$ gets smaller, $\mathcal{E}_l$ is
equipped with the stronger weight. Indeed, $\mathcal{E}_0$ corresponds
to the positive part of the real energy $E$.  First we rewrite the Euler-Posisson
system (1.4) and (1.5) for $\gamma=\frac{6}{5}$ by letting
$\rho=\rho_0+\sigma$ :
\begin{equation}
\sigma_t + \frac{1}{r^2}(r^2\rho_0 u)_r+ \frac{1}{r^2}(r^2\sigma
u)_r =0
\end{equation}
\begin{equation}
u_t+ uu_r+
\frac{4\pi}{15}(\rho_0+\sigma)^{-\frac{4}{5}}(\rho_0+\sigma)_r
+\frac{4\pi}{r^2}\int_0^r(\rho_0+\sigma)s^2 ds=0
\end{equation}
If we define $\mathcal{E}_0^0(t)$ and $T^{\delta}$ by
\begin{equation*}
\begin{split}
\text{ }\text{ }\text{ }\text{ }\text{ }\mathcal{E}_0^0(t)&=\int_{\mathbb{R}^3}(\rho_0+\sigma)u^2+
\frac{4\pi}{15}(\rho_0+\sigma)
^{-\frac{4}{5}}\sigma^2 dx\text{ }\text{ and }\\
\text{ }\text{ }\text{ }\text{ }\text{ }T^{\delta}\text{ }&=
\frac{1}{\sqrt{\mu_0}}\ln \frac{\theta}{\delta} \text{ }\text{ where }
\sqrt{\mu_0} \text{ is the sharp linear growth rate,}
\end{split}
\end{equation*}
the main
theorem can be stated as follows :\\

\begin{thm} Let $l\leq -3$ be fixed. There exist $\theta>0$ and $C>0$ such that for any
small $\delta>0$, there exists a family of solutions
$\nu^{\delta}(t)=\binom{\sigma^{\delta}(t)}{u^{\delta}(t)}$ of
(1.9) and (1.10) such that
\begin{equation}
\sqrt{\mathcal{E}_l}(0)\leq C\delta, \text{  but } \sup_{0\leq t
\leq T^{\delta}}\sqrt{\mathcal{E}_0^0}(t)\geq \theta.
\end{equation}
\end{thm}\

We remark that the escape time $T^{\delta}$ is determined by the
exponential growth rate of the linearized system (2.1) and (2.2)
and the instability occurs before the possible blowup of the
smooth solutions. The local existence of regular solutions for
$\gamma=\frac{6}{5}$ including (1.8) was shown in \cite{ga} and
\cite{Bez} independently. Gamblin uses the paradifferential
calculus and Bezard uses the Gagliardo-Nirenberg inequality and
Littlewood-Paley theory. In this paper we take their existence
results for granted without proving again.

Rein's work in \cite{r} implies that the steady states for
$\gamma>\frac{4}{3}$ are stable under the energy functional $E$
as minimizers. By contrast with those cases, our steady profile
for $\gamma=\frac{6}{5}$ does not minimize $E$ and that
motivated us to investigate an instability under the energy-like
measurement $\mathcal{E}_0^0$. The proof of Theorem 1.1 is based on
the bootstrap argument from the linear instability to the
nonlinear dynamical model. In general, passing a linearized
instability to a nonlinear instability needs much effort in the
PDE context: the spectrum of the linear part is fairly complicated
and the unboundedness of the nonlinear part usually yields a loss
in derivatives. In order to get around these difficulties, the
careful analysis on the linearized system which controls the
sharp exponential growth rate of its solutions is necessary.
In addition, the energy estimates on the whole system and the
interplay with the linear analysis can close the argument. See
\cite{gs1} for its original method.

The main difficulty in this paper is to derive Proposition 6.1, a
key estimate for the bootstrap argument. There are two important
ingredients: the first idea is to find the form of
\begin{equation}\frac{d}{dt}\mathcal{E}\leq \eta\mathcal{E}+
(\mathcal{E})^2 + \text{(quadratic but lower derivative terms) }
\end{equation}
where $\eta$ is smaller than the sharp growth rate
$\sqrt{\mu_0}$, but it turns out that (1.12) is not enough to
close the energy estimate; so we introduce the weighted energy
$\mathcal{E}_l$ by the utilization of a family of symmetrizers as
weights. New undesirable quadratic terms come out during the weighted energy
estimates, but they turn out to have weaker weights. As for weighted
cubic terms, we introduce the weighted Gagliardo-Nirenberg inequality
(6.19). Here is the modified version of (1.12):
\begin{equation}
\begin{split}
\frac{d}{dt}\mathcal{E}_l&\leq \eta\mathcal{E}_l+
(\text{quadratic but better weighted terms})\\&+
(\mathcal{E}_l)^2 + \text{(quadratic but lower derivative terms)}
\end{split}
\end{equation}
The potential part
$\nabla\Phi$ not only has the smoothing effect to the whole system
but also behaves nicely with respect to weights. Eventually, in
cooperating with weights, the chain type estimates by using (1.13)
complete the bootstrap argument
as well as the proof of Proposition 6.1.

Another technical difficulty lies in that both linearized and
full system do not comply submissively with spatial derivatives.
Our argument in both linear and nonlinear parts heavily depends on
the estimates on pure temporal derivative terms. As for spatial and
mixed derivative terms we use the equations directly to control
them in terms of temporal derivative terms. Furthermore, when one
uses the polar coordinates, we seem to end up with simpler one
dimensional flow but the singularity at the origin comes into
play. For clear and precise understanding of the largest growing
mode, we take the polar coordinates in the linear analysis while
the rectangular ones are used for the nonlinear analysis to take
care of spatial and mixed derivative terms.

Other cases of $\gamma$, $\text{ }\frac{6}{5}<\gamma<\frac{4}{3}$
attain more physical, interesting feature involving the vacuum
boundary. For those cases, steady states satisfying (1.6) are
compactly supported and even the local existence of the
Euler-Poisson system including those stationary solutions has not
been completed yet. At this moment the above argument does not
seem to apply directly to such $\gamma$'s. But we believe the
method developed in this paper can make a contribution to show the
instability for those $\gamma$'s. We will leave them for future study.

The paper proceeds as follows. The first half of this article is
devoted to develop the linear theory: finding the largest growing
mode of the linearized Euler-Poisson system and deriving some
regularity of it. While the linear instability was studied in \cite{lin}, it was done in Lagrangian formulation and it does not give
the precise growth rate. For our purpose we shall demonstrate the
explicit linear instability analysis in Eulerian coordinates. In
Section 2, we formulate a variational problem to find the biggest
eigenvalue $\sqrt{\mu_0}$ and corresponding eigenfunction. In the
subsequent section, the fast decay property of the largest
growing mode is derived. In Section 4, we show that
$\sqrt{\mu_0}$ dominates exponential growth rate of any solutions
for the linearized system. In the other half the nonlinear
analysis is carried out. In Section 5, the weighted instant energy
and total energy are introduced and it will be shown that the
total energy is bounded by the instant energy containing only
temporal derivative terms. We perform the weighted energy
estimates in Section 6 to get the precise estimate of (1.13).
Finally in Section 7, the
bootstrap argument and Theorem 1.1 will be proven.\\

\section{Existence of the Largest Growing Mode}

Firstly, we study the linearized Euler-Poisson equations for the
spherically symmetric case. We assume $\rho>0$ which is our
concern. Letting $\rho=\rho_0+\sigma$, we linearize (1.4) and
(1.5) around a given steady state $(\rho_0, 0)$ and get the
linearized Euler-Poisson equations in terms of $\sigma$ and u.

\begin{equation}
\frac{\partial\sigma}{\partial t} + \frac{\partial\rho_0}{\partial
r}u + \rho_0\frac{\partial u}{\partial r} + \frac{2}{r}\rho_0 u =
0
\end{equation}
\begin{equation}
\frac{\partial u}{\partial t} +A \{\gamma\rho_0^{\gamma-2}
\frac{\partial\sigma}{\partial r}
+\gamma(\gamma-2)\rho_0^{\gamma-3} \frac{\partial\rho_0}{\partial
r}\sigma\} + \frac{4\pi}{r^2} \int_0^r \sigma s^2 ds = 0.
\end{equation}\

To find a growing mode for (2.1) and (2.2), let $\sigma=e^{\lambda
t}\phi(r)$ and $u=e^{\lambda t}\psi(r)$. Then (2.1) becomes
\begin{equation}
\lambda\phi + \frac{\partial\rho_0}{\partial r}\psi + \rho_0
\frac{\partial\psi}{\partial r} + \frac{2}{r}\rho_0\psi=0.
\
\end{equation}
Since $\frac{\partial\rho_0}{\partial r}\psi + \rho_0
\frac{\partial\psi}{\partial r} + \frac{2}{r}\rho_0\psi =
\frac{1}{r^2}(r^2\rho_0\psi)_r$, (2.3) gives a simple relation
between $\phi$ and $\psi$:
\begin{equation}
\psi(r)=-\frac{\lambda}{r^2\rho_0}\int_{0}^{r}\phi(s)s^2ds.
\end{equation}
Similarly, (2.2) becomes
\begin{equation}
\lambda\psi + A[\gamma\rho_0^{\gamma-2}\frac{\partial\phi}
{\partial r} + \gamma(\gamma-2)\rho_0^{\gamma-3}\frac{\partial
\rho_0}{\partial r}\phi]+\frac{4\pi}{r^2}\int_{0}^{r}\phi s^2ds
=0.
\end{equation}
Multiplying (2.5) by $\frac{\lambda}{A}$ and using (2.3), we get
\begin{equation*}
\begin{split}
\frac{\lambda^2\psi}{A}=&\gamma\rho_0^{\gamma-2}(-\lambda\phi)'
+\gamma(\gamma-2)\rho_0^{\gamma-3}\rho_0'(-\lambda\phi)+ \frac
{4\pi\rho_0}{A}(-\frac{\lambda}{r^2\rho_0}\int_0^{\infty}\phi(s)ds)\\
=&\gamma\rho_0^{\gamma-2}[\rho_0''\psi+2\rho_0'\psi'+\rho_0\psi''
-\frac{2\rho_0\psi}{r^2}+\frac{2\rho'\psi}{r}+\frac{2\rho_0\psi'}{r}]\\
&+\gamma(\gamma-2)\rho_0^{\gamma-3}\rho_0'[\rho_0'\psi+\rho_0\psi'
+\frac{2\rho_0\psi}{r}]+\frac{4\pi\rho_0}{A}\psi\\
=&\gamma\rho_0^{\gamma-1}\psi''+[\gamma^2\rho_0^{\gamma-2}\rho_0'
+\frac{2\gamma}{r}\rho_0^{\gamma-1}]\psi'+[\gamma\rho_0^{\gamma-2}
\rho_0''-\frac{2\gamma}{r^2}\rho_0^{\gamma-1}+\frac{2\gamma}{r}
\rho_0^{\gamma-2}\rho_0'\\
&+\gamma(\gamma-2)\rho_0^{\gamma-3}
(\rho_0')^2+\frac{2\gamma(\gamma-2)}{r}\rho_0^{\gamma-2}\rho_0'
+\frac{4\pi\rho_0}{A}]\psi
\end{split}
\end{equation*}
where $'=\frac{d}{dr}$ and this is the 2nd order ordinary
differential equation. For the further simplification recall that
$\rho_0$ satisfies (1.6): $A\gamma\rho_0^{\gamma-2}\rho_0'
+\frac{4\pi}{r^2}\int_0^r\rho_0 s^2 ds =0$. Compute
$\frac{1}{Ar^2}\frac{d}{dr}[r^2\cdot(1.6)]$, and then we get the following
relation:
\begin{equation*}
\gamma\rho_0^{\gamma-2}
\rho_0''+\frac{2\gamma}{r}
\rho_0^{\gamma-2}\rho_0'+\gamma(\gamma-2)\rho_0^{\gamma-3}
(\rho_0')^2+\frac{4\pi\rho_0}{A}=0.\\
\
\end{equation*}
In turn we have
\begin{equation}
\frac{\lambda^2\psi}{A}=\gamma\rho_0^{\gamma-1}\psi'' +
[\gamma^2\rho_0^{\gamma-2}\rho_0' +
\frac{2\gamma}{r}\rho_0^{\gamma-1}]\psi' +
[\frac{2\gamma(\gamma-2)}{r}\rho_0^{\gamma-2}\rho_0' -
\frac{2\gamma}{r^2}\rho_0^{\gamma-1}]\psi.
\end{equation}
Multiply (2.6) by $\frac{\rho_0 r^2}{\gamma}$ to obtain the
following:

\begin{equation}
\begin{split}
\lambda^2\frac{\rho_0 r^2}{A\gamma}\psi&=r^2 \rho_0^{\gamma}\psi''
+[\gamma r^2 \rho_0^{\gamma-1}\rho_0' + 2r\rho_0^{\gamma}]\psi' +
[2(\gamma-2)r\rho_0^{\gamma-1}\rho_0'-2\rho_0^{\gamma}]\psi\\
 &= r^2 \rho_0^{\gamma}\psi'' + (r^2\rho_0^{\gamma})'\psi' +
2[\frac{\gamma-2}{\gamma}r(\rho_0^{\gamma})'-\rho_0^{\gamma}]\psi\\
&= (r^2\rho_0^{\gamma}\psi')'+ 2[\frac{\gamma-2}{\gamma}r
(\rho_0^{\gamma})'-\rho_0^{\gamma}]\psi\\
\
\end{split}
\end{equation}
Denote the RHS of (2.7) by $L\psi$:
\begin{equation}
L\psi\equiv(r^2\rho_0^{\gamma}\psi')'+ 2[\frac{\gamma-2}{\gamma}r
(\rho_0^{\gamma})'-\rho_0^{\gamma}]\psi
\end{equation}
Note that the linear operator $L$ is self-adjoint and hence
$\lambda^2$ is real.\\

\begin{lem}
Suppose $\chi$ and $\omega$ ($\omega>0$) satisfy
$L\chi=\omega^2\frac{\rho_0 r^2}{\gamma}\chi$. Define
$\varphi=-\frac{1}{\omega r^2}(r^2 \rho_0\chi)'$. We assume
$\chi$ and $\varphi$ are well defined admissible functions in a
suitable sense which will be clarified later on. Then
$\sigma=e^{\pm\omega t}\varphi$ and $u=e^{\pm\omega t}\chi$ are a
solution pair of the linearized equations (2.1) and (2.2).
\end{lem}

Proof. This is obvious by the definition of $\chi$ and $\varphi$.
One can keep track of the derivation of $L$ to see it.$\square$\\

Lemma 2.1. tells us that $(e^{\omega t}\varphi,e^{\omega t} \chi)$
satisfying all the assumptions is a growing mode for the
linearized Euler-Poisson equations. Next we show such a growing
mode actually exists when $\gamma=\frac{6}{5}$. This can be done by
looking at the eigenvalue problem of the operator $L$ due to
Lemma 2.1. In other words, we only need to find $\psi$ and
$\lambda$ ($\lambda>0$) such that $L\psi=\lambda^2 \frac{\rho_0
r^2}{A\gamma}\psi$, where $L\psi$ is defined in (2.8).\\

From now on we fix $\gamma=\frac{6}{5}$ and $A=\frac{2\pi}{9}$ and
corresponding $$\rho_0(r)=\frac{1}{(1+r^2)^{\frac{5}{2}}}.$$

The starting equation is $L\psi=\mu \frac{\rho_0 r^2}{A\gamma}\psi$ on
$(0,\infty)$, where $\mu=\lambda^2$. It is well known that the
largest eigenvalue $\mu_0$ is given by a variational formula:

$$\mu_0=\sup [\frac{Q(\psi)}{I(\psi)}: Q(\psi)<\infty,
I(\psi)<\infty]$$ where
\begin{equation*}
\begin{split}
Q(\psi)=(L\psi,\psi)&=-\int_0^{\infty} r^2
\rho_0^{\gamma}(\psi')^2 dr +
2\int_0^{\infty}(\frac{\gamma-2}{\gamma}r
(\rho_0^{\gamma})'-\rho_0^{\gamma})\psi^2 dr\\
&=-\int_0^{\infty}\frac{r^2}{(1+r^2)^3}(\psi')^2 dr +
2\int_0^{\infty}\frac{3r^2-1}{(1+r^2)^4}(\psi)^2 dr\\
\end{split}
\end{equation*}
and
$$I(\psi)=(\frac{\rho_0 r^2}{A\gamma}\psi, \psi)=\int_0^{\infty}
\frac{\rho_0 r^2}{A\gamma}\psi^2 dr
=\frac{15}{4\pi}\int_0^{\infty}\frac{r^2}{(1+r^2)^{\frac{5}{2}}}(\psi)^2
dr.$$\

Hence once the above formula attains the $\sup$, the
largest eigenvalue $\mu_0$ and corresponding eigenfunction
$\psi_0$ of the linear operator $L$ gives a largest growing
mode of the linearized equations. In order to carry it out, first
define a norm for any $\psi\in\mathcal{C}_c^{\infty}(0,\infty),$
$$\|\psi\|^2 \equiv
\int_0^{\infty}\frac{r^2}{(1+r^2)^3}(\psi')^2 dr +
\int_0^{\infty} \frac{2}{(1+r^2)^4}\psi^2 dr + \frac{15}{4\pi}
\int_0^{\infty}\frac{r^2}{(1+r^2)^{\frac{5}{2}}}\psi^2 dr. $$
Let $\mathrm{H}=\overline{\mathcal{C}_c^{\infty}(0,\infty)}$ in
the above norm. It is clear that $\frac{Q(\psi)}{I(\psi)}=
\frac{Q(c\psi)} {I(c\psi)}$ for any nonzero constant $c$. Thus
the variatonal problem can be rephrased as to find a maximum
$\mu_0$ of $Q(\psi)$ on $\mathrm{H}$ under the normalization
condition $I(\psi)=1$.\\

\begin{prop} There exists a $\psi_0\in\mathrm{H}$ such that
$I(\psi_0)=1$ and $Q(\psi_0)=\mu_0$ ($\mu_0>0$) i.e. the $\sup$
$\mu_0$ is attained on $\mathrm{H}$.
\end{prop}

Proof: First we claim $\mu_0>0$. Consider $\psi=\sqrt{r}$.

$$\frac{Q(\sqrt{r})} {I(\sqrt{r})} =
\frac{-\int_0^{\infty}\frac{r^2}{(1+r^2)^3}(\frac{1}{2\sqrt{r}})^2
dr + 2\int_0^{\infty}\frac{3r^2-1}{(1+r^2)^4}(\sqrt{r})^2
dr}{\frac{15}{4\pi}\int_0^{\infty}\frac{r^2}{(1+r^2)^{\frac{5}{2}}}(\sqrt{r})^2
dr}$$

Since $\mu_0\geq\frac{Q(\sqrt{r})} {I(\sqrt{r})}$ and
$I(\sqrt{r})>0$, it is enough to show $Q(\sqrt{r})>0$.
\begin{equation*}
\begin{split}
Q(\sqrt{r})&=(-\frac{1}{4}+6)\int_0^{\infty}\frac{r}{(1+r^2)^3} dr
- 8\int_0^{\infty}\frac{r}{(1+r^2)^4}dr\\
&=\frac{23}{4}[-\frac{1}{2(1+r^2)^2}]_0^{\infty}-8
[-\frac{1}{3(1+r^2)^3}]_0^{\infty}\\
&=\frac{23}{8}-\frac{8}{3} >0\\
\end{split}
\end{equation*}\

And note that the positive part of $Q$ is uniformly bounded by
$I$ because $\frac{3r^2}{(1+r^2)^4}=O(\frac{1}{r^6})$ and $\frac
{r^2}{(1+r^2)^{\frac{2}{5}}}=O(\frac{1}{r^3})$ for sufficiently
large $r$. This implies $\mu_0$ is finite. To show $\mu_0$ is
attained on $\mathrm{H}$ let $\{\psi_n\}$ be a maximizing
sequence i.e. $$Q(\psi_n)\nearrow\mu_0 \text{  as }
n\longrightarrow\infty \text{ and } I(\psi_n)=1 \text{ for all }
n.$$ Let $\psi_0$ be its weak limit. Then by the lower
semicontinuity of weak convergence, we have $$\liminf
\int\frac{r^2}{(1+r^2)^3}(\psi_n')^2 dr \geq
\int\frac{r^2}{(1+r^2)^3}(\psi_0')^2 dr,$$
$$\liminf\int\frac{1}{(1+r^2)^4}\psi_n^2 dr \geq
\int\frac{1}{(1+r^2)^4}\psi_0^2 dr,$$
$$\liminf \int\frac{r^2}{(1+r^2)^{\frac{5}{2}}}\psi_n^2 dr\geq
\int\frac{r^2}{(1+r^2)^{\frac{5}{2}}}\psi_0^2 dr.$$\\

Claim 1. (Compactness of the positive part) There exists a
subsequence $\{\psi_{n_k}\}$ of $\{\psi_n\}$ such that
$$\int\frac{r^2}{(1+r^2)^4}\psi_{n_k}^2 dr \longrightarrow
\int\frac{r^2}{(1+r^2)^4}\psi_0^2 dr \text{ as }
n_k\longrightarrow\infty.$$\

Claim 2. $I(\psi_0)=1$. \\

Since Claim 1 immediately implies $Q(\psi_0)=\mu_0$, the
conclusion follows from Claim 1 and Claim 2. It remains to prove
Claim 1 and Claim 2. Claim 2 follows from a simple scaling
argument. Suppose
$\frac{15}{4\pi}\int\frac{r^2}{(1+r^2)^{\frac{5}{2}}}\psi_0^2
dr=\alpha^2<1$. Then $$I(\frac{\psi_0}{\alpha})=1 \text{ and }
Q(\frac{\psi_0}{\alpha})=\frac{1}{\alpha^2}\mu_0>\mu_0$$ which is
a contradiction to the definition of $\mu_0$. To prove Claim 1,
first observe that $$\int_R^{\infty}\frac{r^2}{(1+r^2)^4}\psi_n^2
dr\leq\frac{1}{(1+R^2)^{\frac{3}{2}}}\int_R^{\infty}
\frac{r^2}{(1+r^2)^{\frac{5}{2}}}\psi_n^2
dr\leq\frac{4\pi}{15(1+R^2)^{\frac{3}{2}}}.$$
Fix $R>0$. On the finite interval $(0,R)$, since $\int_0^R
\frac{r^2}{(1+r^2)^4}\psi^2 dr \sim \|\psi\|_{L^2(B_R(0))}^2$,
$\int_0^R \frac{r^2}{(1+r^2)^\frac{5}{2}}\psi^2 dr \sim
\|\psi\|_{L^2(B_R(0))}^2$, and $\int_0^R
\frac{r^2}{(1+r^2)^3}(\psi')^2 dr \sim \|\psi'\|_{L^2(B_R(0))}^2$
where $B_R(0)$ is a ball with radius $R$ in $\mathbb{R}^3$, we can
apply the Rellich-Kondrachov Compactness theorem, which says
$H^1(B_R(0))$ is compactly embedded in $L^q(B_R(0))$ for each
$1\leq q<6$, . So there exists a subsequence $\{\psi_{n_k}\}$
such that

\begin{equation*}
\int_0^R\frac{r^2}{(1+r^2)^4}\psi_{n_k}^2 dr\longrightarrow
\int_0^R\frac{r^2}{(1+r^2)^4}\psi_0^2 dr, \text{ as }
n_k\longrightarrow\infty.
\end{equation*}
For any given $\epsilon>0$, choose $R>0$ large enough so that
$\frac{4\pi}{15(1+R^2)^{\frac{3}{2}}}< \epsilon$.

\begin{equation*}
\begin{split}
\int_0^{\infty}\frac{r^2}{(1+r^2)^4}\psi_{n_k}^2 dr &= \int_0^R
\frac{r^2}{(1+r^2)^4}\psi_{n_k}^2 dr
+\int_R^{\infty}\frac{r^2}{(1+r^2)^4}\psi_{n_k}^2 dr \\ &<\int_0^R
\frac{r^2}{(1+r^2)^4}\psi_{n_k}^2 dr +\epsilon
\end{split}
\end{equation*}
Now take the limit of $n_k\longrightarrow\infty$. Since
$\epsilon>0$ is arbitrary, this finishes the proof of Proposition
2.1.$\square$\\

Finally, let us make sure $\psi_0$ is in fact an eigenfunction
corresponding to $\mu_0$. Consider a perturbation $\psi_0 +
\epsilon\eta$ around $\psi_0$. Since $\mu_0=Q(\psi_0)$, $$\frac{Q(
\psi_0 + \epsilon\eta)}{I(\psi_0 + \epsilon\eta)}\leq
Q(\psi_0)=\mu_0$$ for all sufficiently small $\epsilon$ and all
admissible function $\eta$. And  hence we have

\begin{equation*}
\begin{split}
Q(\psi_0 + \epsilon\eta)&\leq Q(\psi_0)I(\psi_0 + \epsilon\eta)\\
Q(\psi_0)+2\epsilon(L\psi_0,\eta)+\epsilon^2Q(\eta)&\leq
Q(\psi_0)\{I(\psi_0)+2\epsilon(\frac{\rho_0
r^2}{\gamma}\psi_0,\eta)+\epsilon^2 I(\eta)\}
\
\end{split}
\end{equation*}
Here $(\cdot,\cdot)$ is the standard inner product in
$L^2(0,\infty)$. Note $I(\psi_0)=1$ and $Q(\eta)\leq\mu_0
I(\eta)$. Hence in order for the above inequality to hold for all
$\epsilon$ and $\eta$, the coefficient of $\epsilon$ should
vanish, i.e. $(L\psi_0, \eta)=\mu_0(\frac{\rho_0
r^2}{A\gamma}\psi_0,\eta)$ for all $\eta$. Thus,
$L\psi_0=\mu_0\frac{\rho_0r^2}{A\gamma}\psi_0$ and therefore
$\sqrt{\mu_0}$ and $\psi_0$ give a
largest growing mode for the linearized Euler-Poisson equations.\\

\section{The Regularity of the Largest Growing Mode}

Recall $\psi_0$ satisfies the following 2nd linear ordinary
differential equation:

\begin{equation}
\begin{split}
\frac{15\mu_0}{4\pi}\frac{r^2}{(1+r^2)^{\frac{5}{2}}}\psi_0 &=
(\frac{r^2}{(1+r^2)^3}\psi_0')' +
2\frac{3r^2-1}{(1+r^2)^4}\psi_0\\
&=\frac{r^2}{(1+r^2)^3}\psi_0'' + \frac{2r-4r^3}{(1+r^2)^4}\psi_0'
+2\frac{3r^2-1}{(1+r^2)^4}\psi_0\\
\
\end{split}
\end{equation}

The interior regularity easily follows from the elliptic theory
of the 2nd order differential equations: since each coefficient
of $\psi_0''$, $\psi_0'$ and $\psi_0$ in (3.1) is in
$\mathrm{C}^{\infty}(\epsilon, R)$ where $\epsilon$ is a small
enough positive fixed number and $R$ is a large enough fixed
number, $\psi_0$ is also $\mathrm{C}^{\infty}$ on $(\epsilon,
R)$. In the following two subsections we investigate the behavior
of $\psi_0$ when $r$ is either very small or very large.

\subsection{The behavior of $\psi_0$ near the origin}

In this first subsection we show $\psi_0$ is analytic near the
origin. This rather surprising property easily follows from the
classical theorem by Frobenius from the ODE theory. For
reference, see \cite{braun}. Before we prove the analyticity, it will
be shown that the maximizing property implies
the boundedness of $\psi_0$ at the origin.\\

\begin{lem} There exists a decreasing sequence $\{\epsilon_k\}$
with $\epsilon_k\searrow 0$ such that
$\psi_0(\epsilon_k)\psi_0'(\epsilon_k)\geq 0$ for each k.
\end{lem}

Proof. Consider $\psi_1=\Theta\psi_0$, where $\Theta$ is a
lipschitz cutoff function defined by
\[
\Theta(r)=
  \begin{cases}
    0 & \text{if $0\leq r \leq \epsilon$}, \\
    r-\epsilon & \text{if $\epsilon\leq r \leq 2\epsilon$},\\
    \epsilon & \text{if $r\geq\ 2\epsilon$}.\\
  \end{cases}
\]
Here, $\epsilon$ is a sufficiently small positive number to be
clarified. Then clearly $\psi_1$ is in the admissible set
$\mathrm{H}$. First note that $I(\psi_1)\leq I(\epsilon \psi_0)$
since $\psi_1^2 \leq \epsilon^2\psi_0^2$. Next let us look at the
difference between $Q(\psi_1)$ and $Q(\epsilon\psi_0)$. Since
$\psi_1$ has the same value as $\epsilon\psi_0$ on $r\geq
2\epsilon$ and
$\psi_1'^2=\psi_0^2\mathbf{1}_{(\epsilon,2\epsilon)}+
2(r-\epsilon)\psi_0\psi_0'\mathbf{1}_{(\epsilon,2\epsilon)} +
(r-\epsilon)^2\psi_0'^2\mathbf{1}_{(\epsilon,2\epsilon)} +
\epsilon^2\psi_0'^2\mathbf{1}_{(2\epsilon,\infty)}$,
$Q(\psi_1)-Q(\epsilon \psi_0)$ is
\begin{equation}
\begin{split}
& -\int_{\epsilon}^{2\epsilon} \frac{r^2}{(1+r^2)^3}[\psi_0^2 +
(r-\epsilon)^2\psi_0'^2 + 2(r-\epsilon)\psi_0\psi_0']dr +
2\int_{\epsilon}^{2\epsilon}
\frac{3r^2-1}{(1+r^2)^4}(r-\epsilon)^2 \psi_0^2 dr \\ &+
\int_0^{2\epsilon}\frac{r^2}{(1+r^2)^3}\epsilon^2 \psi_0'^2 dr -2
\int_0^{2\epsilon}\frac{3r^2-1}{(1+r^2)^4}\epsilon^2\psi_0^2 dr.
\
\end{split}
\end{equation}
After collecting similar terms, (3.2) can be written as following:
\begin{equation}
\begin{split} &
\int_0^{\epsilon}\frac{r^2}{(1+r^2)^3}\epsilon^2\psi_0'^2 dr +
\int_{\epsilon}^{2\epsilon}\frac{r^2}{(1+r^2)^3}r(2\epsilon-r)\psi_0'^2
dr
-\int_{\epsilon}^{2\epsilon}\frac{2r^2(r-\epsilon)}{(1+r^2)^3}\psi_0\psi_0'dr \\
&+ \int_0^{\epsilon}\frac{2-6r^2}{(1+r^2)^4}\epsilon^2\psi_0^2 dr
+ \int_{\epsilon}^{2\epsilon}\frac{5r^4-12\epsilon r^3 -3r^2 +
4r\epsilon}{(1+r^2)^4}\psi_0^2 dr.
\
\end{split}
\end{equation}
Call the first three terms of (3.3) (I) and the last two terms
(II). Now suppose the lemma were false. Then for any small enough
$\epsilon>0$, we may assume $\psi_0\psi_0' < 0$ on
$(0,2\epsilon)$. It is obvious to see that (I) is positive. And
because $\psi_0^2$ is decreasing on $(0,2\epsilon)$,
$$\int_0^{\epsilon}\frac{2-6r^2}{(1+r^2)^4}\epsilon^2\psi_0^2 dr >
\int_{\epsilon}^{2\epsilon}\frac{2-6r^2}{(1+r^2)^4}\epsilon^2\psi_0^2
dr $$ and hence
\begin{equation*}
\begin{split}
\text{(II)}&>\int_{\epsilon}^{2\epsilon}\frac{2-6r^2}{(1+r^2)^4}
\epsilon^2\psi_0^2 dr +
\int_{\epsilon}^{2\epsilon}\frac{5r^4-12\epsilon r^3 -3r^2 +
4r\epsilon}{(1+r^2)^4}\psi_0^2 dr\\
&=\int_{\epsilon}^{2\epsilon}\frac{5r^4-12\epsilon
r^3-6\epsilon^2 r^2+\epsilon^2}{(1+r^2)^4}\psi_0^2 dr +
\int_{\epsilon}^{2\epsilon}\frac{\epsilon^2+4\epsilon r
-3r^2}{(1+r^2)^4}\psi_0^2 dr\\
&=\text{(III)} +\text{(IV)}.\\
\end{split}
\end{equation*}

(III) $\geq 0$ since $5r^4-12\epsilon r^3-6\epsilon^2
r^2+\epsilon^2\geq 0$ on $(\epsilon, 2\epsilon)$ for sufficiently
small $\epsilon >0$.\\

Claim. (IV)$\geq 0$.\\

Claim implies $Q(\psi_1)-Q(\epsilon\psi_0)=$(I)+(II)$> 0$. This
leads to
$$\frac{Q(\psi_1)}{I(\psi_1)} > \frac{Q(\epsilon\psi_0)}
{I(\epsilon\psi_0)}=\frac{Q(\psi_0)}{I(\psi_0)}=\mu_0.$$ which is
a contradiction to the definition of $\mu_0$. Therefore, to finish
the proof it suffices to verify Claim. Let $\alpha$ be a zero of
$\epsilon^2+4\epsilon r -3r^2=0$ with
$\epsilon<\alpha<2\epsilon$. Then
\[
\epsilon^2+4\epsilon r -3r^2=
  \begin{cases}
    \geq 0& \text{if $\epsilon\leq r\leq\alpha$}, \\
    \leq 0& \text{if $\alpha \leq r \leq 2\epsilon$}.\\
  \end{cases}
\]
The decreasing assumption on $\psi_0^2$ will be again used:
\begin{equation*}
\begin{split}
\text{(IV)}&\geq \min_{[\epsilon, \alpha]}
\frac{\psi_0^2(r)}{(1+r^2)^4}\int_{\epsilon}^{\alpha}
\epsilon^2+4\epsilon r -3r^2 dr + \max_{[\alpha, 2\epsilon]}
\frac{\psi_0^2(r)}{(1+r^2)^4}\int_{\alpha}^{2\epsilon}
\epsilon^2+4\epsilon r -3r^2 dr \\
&=\frac{\psi_0^2(\alpha)}{(1+\alpha^2)^4}\int_{\epsilon}^{2\epsilon}
\epsilon^2+4\epsilon r -3r^2 dr\\
&=0.\\
\end{split}
\end{equation*}
The first equality holds because $\frac{\psi_0^2}{(1+r^2)^4}$ is
decreasing and the second one is simply due to the fact
$\int_{\epsilon}^{2\epsilon} \epsilon^2+4\epsilon r -3r^2 dr=0.$
This completes the argument.$\square$\\

The Frobenius theorem with Lemma 3.1 gives rise to the following
analytic property of $\psi_0$ around the origin.

\begin{lem}$\psi_0$ is analytic at $r=0$ and moreover $\psi_0(r)=
ar+o(r^2)$ around the origin where $a$ is a constant.
\end{lem}

Proof. In order to employ the Frobenius theorem, we need to show
the equation (3.1) has a regular singular point at $r=0$. To check
this out in the context of \cite{braun} (p. 215), we rewrite (3.1) in the
following form
\begin{equation*}
\psi_0''+ \frac{2-4r^2}{r(1+r^2)}\psi_0'+
\{\frac{2}{r^2}\frac{3r^2-1}{1+r^2}-\frac{15\mu_0}{4\pi}(1+r^2)
^{\frac{1}{2}}\}\psi_0=0.
\end{equation*}

Let $P(r)$ and $Q(r)$ be coefficients of $\psi_0'$ and $\psi_0$
respectively in the above. Then it is clear that $rP(r)$ and
$r^2Q(r)$ are analytic at $r=0$, which means $r=0$ is a regular
singular point of (3.1). Let $p_0$, $q_0$ be the zeroth order term
of $rP(r)$ and $r^2Q(r)$ respectively. It is easy to check
$p_0=2$ and $q_0=-2$. The indicial equation $r(r-1)+p_0 r+q_0=0$
has two roots $r_1=1$ or $r_2=-2$. Hence by the Frobenius theorem
$\psi_0$ has a power series representation of either
$y_1(r)=r\sum_{n=0}^{\infty}a_n r^n$ or $y_2(r)=a y_1(r) \ln r+
r^{-2}\sum_{n=0}^{\infty}b_n r^n$. However, $y_2(r)$ is impossible
by Lemma 3.1 and therefore
$\psi_0=r\sum_{n=0}^{\infty}a_n r^n$ is obtained.$\square$

\subsection{Asymptotic behavior of $\psi_0$}

In this subsection we first observe $\psi_0^2(r)$ is nonincreasing
near the infinity owing to the maximizing property and then the
fast decay of $\psi_0$ in an appropriate sense will be shown by
using a standard bootstrap argument.\\

\begin{lem}
There exists an increasing sequence $\{R_k\}$ with
$R_k\nearrow\infty$ such that $\psi_0(R_k)\psi_0'(R_k)\leq
0$ for each $k$.
\end{lem}

Proof. Consider $\psi_1 =\Theta\psi_0$, where $\Theta$ is a
lipschitz cutoff function defined by
\[
\Theta(r)=
  \begin{cases}
    1 & \text{if $0\leq r \leq R$}, \\
    R+1-r & \text{if $R\leq r \leq R+1$},\\
    0 & \text{if $r\geq\ R+1$}.\\
  \end{cases}
\]
Here, $R$ is a sufficiently large number to be determined. Note
that $(\psi_1')^2=(\psi_0')^2\mathbf{1}_{(0,R)}+(R+1-r)^2
(\psi_0')^2\mathbf{1}_{(R, R+1)}+\psi_0^2\mathbf{1}_{(R,
R+1)}-2(R+1-r)\psi_0\psi_0'\mathbf{1}_{(R, R+1)}$. Then clearly
$\psi_1\in\mathrm{H}$. Let us compute $I(\psi_1)$ and $Q(\psi_1)$.
Recall $I(\psi_0)=1$ and $Q(\psi_0)=\mu_0$.
\begin{equation*}
\begin{split}
 I(\psi_1)=&I(\psi_0)-\frac{15}{4\pi}\int_R^{\infty}\frac{r^2}
 {(1+r^2)^{\frac{5}{2}}}\psi_0^2
 dr +\frac{15}{4\pi}
 \int_R^{R+1}(R+1-r)^2\frac{r^2}{(1+r^2)^{\frac{5}{2}}}\psi_0^2
 dr\\
 =&1-\frac{15}{4\pi}[\int_R^{R+1}(1-(R+1-r)^2)\frac{r^2}{(1+r^2)
 ^{\frac{5}{2}}}\psi_0^2
 dr + \int_{R+1}^{\infty}\frac{r^2}{(1+r^2)^{\frac{5}{2}}}\psi_0^2
 dr ]\\
 \equiv& 1-C_R\\
\
Q(\psi_1)=&Q(\psi_0)-[-\int_R^{\infty}\frac{r^2}{(1+r^2)^3}\psi_0'^2
dr + 2\int_R^{\infty}\frac{3r^2-1}{(1+r^2)^4}\psi_0^2 dr \\
&+\int_R^{R+1}\frac{r^2}{(1+r^2)^3}\psi_0^2 dr +
\int_R^{R+1}(R+1-r)^2\frac{r^2}{(1+r^2)^3}\psi_0'^2 dr \\
-&2\int_R^{R+1}(R+1-r)\frac{r^2}{(1+r^2)^3}\psi_0\psi_0' dr -
2\int_R^{R+1}(R+1-r)^2\frac{3r^2-1}{(1+r^2)^4}\psi_0^2 dr]\\
\geq&\mu_0-[-2\int_R^{R+1}(R+1-r)\frac{r^2}{(1+r^2)^3}\psi_0\psi_0'
dr + \int_R^{R+1}\frac{r^2}{(1+r^2)^3}\psi_0^2 dr\\
&+2\int_R^{R+1}(1-(R+s-r)^2)\frac{3r^2-1}{(1+r^2)^4}\psi_0^2 dr +
2\int_{R+1}^{\infty}\frac{3r^2-1}{(1+r^2)^4}\psi_0^2 dr]\\
\equiv &\mu_0 +
2\int_R^{R+1}(R+1-r)\frac{r^2}{(1+r^2)^3}\psi_0\psi_0' dr - D_R\\
\end{split}
\end{equation*}\

Now suppose the proposition were false. Then there exists large
enough $R_0$ such that  if $R>R_0$, then $\psi_0\psi_0'>0$ on
$(R, R+1)$. On the other hand, $D_R\leq\frac{K}{R}C_R$ for
$R>R_0>0$, where $K$ is a fixed constant. So we get the following
inequalities.
\begin{equation*}
\begin{split}
\frac{Q(\psi_1)}{I(\psi_1)}&\geq\frac{\mu_0 +
2\int_R^{R+1}(R+1-r)\frac{r^2}{(1+r^2)^3}\psi_0\psi_0' dr -
D_R}{1-C_R}\\
&> \frac{\mu_0-\frac{K}{R}C_R}{1-C_R}
=\mu_0 + \frac{(\mu_0-\frac{K}{R})C_R}{1-C_R}
\
\end{split}
\end{equation*}
Choosing sufficiently large $R>R_0>0$ so that
$\mu_0\geq\frac{K}{R}$, we get $\frac{Q(\psi_1)}{I(\psi_1)}>\mu_0$
which contradicts the
definition of $\mu_0$. Therefore the lemma follows.$\square$\\

Lemma 3.3 with the integration by parts leads to the next
proposition that shows the fast decay of $\psi_0$ in a suitable
sense. It will play a key role when we prove the main theorem in
Section 7 in a sense that the initial data with different
weights are not essentially different.\\

\begin{prop}
$\int_0^{\infty}\frac{r^{n+2}}{(1+r^2)^{\frac{5}{2}}}\psi_0^2 dr$
and $\int_0^{\infty}\frac{r^{n+2}}{(1+r^2)^3}\psi_0'^2 dr$ are
bounded for each nonnegative integer $n$.
\end{prop}

Proof. Let $R_k$ be given as in Lemma 3.3. Multiply (3.1) by
$r^n\psi_0$ and integrate over $(0,R_k)$. Then we get
\begin{equation}
\frac{15\mu_0}{4\pi}\int_0^{R_k}
\frac{r^{n+2}}{(1+r^2)^{\frac{5}{2}}}\psi_0^2 dr = \int_0^{R_k}r^n
(\frac{r^2}{(1+r^2)^3}\psi_0')'\psi_0 dr +2\int_0^{R_k}
r^n\frac{3r^2-1}{(1+r^2)^4}\psi_0^2 dr
\end{equation}

The first term on the RHS of (3.4) can be rewritten as following
using the integration by parts,
\begin{equation*}
\begin{split}
[&\frac{r^{n+2}}{(1+r^2)^3}\psi_0'\psi_0]_0^{R_k} -
\int_0^{R_k}(nr^{n-1}\psi_0 +
r^n\psi_0')\frac{r^2}{(1+r^2)^3}\psi_0' dr \\
=&\frac{R_k^{n+2}}{(1+R_k^2)^3}\psi_0'(R_k)\psi_0(R_k) -
\int_0^{R_k}
\frac{r^{n+1}}{(1+r^2)^3}[\frac{n}{2}(\psi_0^2)'+(\psi_0')^2] dr\\
=&\frac{R_k^{n+2}}{(1+R_k^2)^3}\psi_0'(R_k)\psi_0(R_k)-\frac{n}{2}
\frac{R_k^{n+1}}{(1+R_k^2)^3}\psi_0^2(R_k)-\int_0^{R_k}
\frac{r^{n+2}}{(1+r^2)^3}(\psi_0')^2 dr \\ &
+\frac{n}{2}\int_0^{R_k}
\frac{(n-5)r^{n+2}+(n+1)r^n}{(1+r^2)^4}\psi_0^2 dr  \\
\end{split}
\end{equation*}
Plugging this into (3.4), we get the following relation.
\begin{equation}
\begin{split}
\frac{15\mu_0}{4\pi}\int_0^{R_k}
\frac{r^{n+2}}{(1+r^2)^{\frac{5}{2}}}\psi_0^2 dr + \int_0^{R_k}
\frac{r^{n+2}}{(1+r^2)^3}&(\psi_0')^2 dr\\
+\frac{R_k^{n+1}}{(1+R_k^2)^3}[\frac{n}{2}\psi_0^2(R_k)-R_k&\psi_0'
(R_k)\psi_0(R_k)]\\
 =\frac{n}{2}\int_0^{R_k}
\frac{(n-5)r^{n+2}+(n+1)r^n}{(1+r^2)^4}\psi_0^2 dr +2&\int_0^{R_k}
r^n\frac{3r^2-1}{(1+r^2)^4}\psi_0^2 dr\
\
\end{split}
\end{equation}\

Observe that for each $R_k$, each term of the LHS of (3.5) is
nonnegative because of the previous lemma and the small $r$ parts
of the RHS of (3.5) are finite due to the behavior of $\psi_0$ near
the origin. First, when $n=0,1,2,3$, all terms in the RHS of (3.5)
are finite for any $R_k$ since they are uniformly bounded by
$\frac{15}{4\pi}\int_0^{\infty}\frac{r^2}{(1+r^2)^{\frac{5}{2}}}\psi_0^2 dr=1$.
As taking the limit of $R_k$ , we know each term of the LHS converges
as $R_k\longrightarrow\infty$. In particular, $\int_0^{\infty}
\frac{r^{n+2}}{(1+r^2)^{\frac{5}{2}}}\psi_0^2 dr$ and $
\int_0^{\infty} \frac{r^{n+2}}{(1+r^2)^3}(\psi_0')^2 dr$ are
bounded for $n=0,1,2,3$. The standard induction on $n$ with (3.4)
gives the desired result for general $n$.$\square$\\

\section{The Linear Growth Rate}

In this section, we show that $\sqrt{\mu_0}$ is the dominating
exponential growth for the linearized Euler-Poisson equations. The
linearized Euler-Poisson equations (2.1) and (2.2) for
$\gamma=\frac{6}{5}$ are
\begin{equation}
\sigma_t + \frac{1}{r^2}(r^2\rho_0 u)_r=0
\end{equation}
\begin{equation}
u_t+\frac{4\pi}{15}(\rho_0^{-\frac{4}{5}}\sigma)_r+\frac{4\pi}{r^2}\int_0^r\sigma
s^2 ds =0\\
\
\end{equation}
Multiply (4.2) by $\frac{15}{4\pi}r^2\rho_0$, take the $t$
derivative and use (4.1) to get rid of $\sigma$.
\begin{equation*}
\begin{split}
\frac{15}{4\pi}r^2\rho_0 u_{tt}=&-r^2\rho_0
(\rho_0^{-\frac{4}{5}}\sigma_t)_r-15\rho_0\int_0^r\sigma_t s^2ds\\
=&r^2\rho_0 (\frac{1}{r^2}\rho_0^{-\frac{4}{5}}(r^2\rho_0
u)_r)_r+15r^2\rho_0^2 u\\
=&r^2\rho_0^{\frac{6}{5}}u_{rr}+\{r^2\rho_0(\rho_0^{\frac{1}{5}})_r
+\rho_0^{\frac{1}{5}}(r^2\rho_0)_r\}u_r\\
&+\{r^2\rho_0(\frac{1}{r^2}
\rho_0^{-\frac{4}{5}}(r^2\rho_0)_r)_r+15r^2\rho_0^2\}u
\end{split}
\end{equation*}
After putting back $\rho_0=\frac{1}{(1+r^2)^{\frac{5}{2}}}$ in
the above, we get an equivalent 2nd order equation for $u\equiv
\Psi$ :
\begin{equation}
\frac{15}{4\pi}\frac{r^2}{(1+r^2)^{\frac{5}{2}}}\Psi_{tt} =
(\frac{r^2}{(1+r^2)^3}\Psi_r)_r + 2(\frac{3r^2-1}{(1+r^2)^4})\Psi
\end{equation}
Here we use $\Psi$ instead of $u$ in order to distinguish the
linear analysis from the nonlinear one. In this way of writing
it is also easy to compare (4.3) with (2.8). Note that $\Psi$ in
(4.3) is a function of both $t$ and $r$ while $\psi$ in (2.8) is a
function of only $r$. Denote the RHS of (4.3) by $L_0\Psi$.
$L_0$ is basically same as $L$ in Section 2.

Define the following quantities.
\begin{equation*}
\begin{split}
&W_0(r)\equiv\frac{15}{4\pi}r^2\rho_0
=\frac{15}{4\pi}\frac{r^2}{(1+r^2)^{\frac{5}{2}}}
\text{  }\\
&\|f\|_W^2\equiv(W f,f)=\int_0^{\infty}Wf^2dr \text{ where } W
\text{ is a given weight function}\\
&P_f(t)\equiv \int_0^{\infty}\frac{r^2}{(1+r^2)^3} f_r^2(t) dr + 2
\int_0^{\infty}\frac{1}{(1+r^2)^4}f^2(t)dr
\end{split}
\end{equation*}\

\begin{lem} For every solution $\Psi$ to (4.3) there exists
$C_{\mu_0}, C_{\mu_0,\alpha}>0$ such that\\ (1)
$$\|\Psi(t)\|_{W_0}\text{ , } \|\Psi_t(t)\|_{W_0} \leq C_{\mu_0}
e^{\sqrt{\mu_0}
t}(\|\Psi_t(0)\|_{W_0}+\|\Psi(0)\|_{W_0}+\sqrt{P_{\Psi}(0)}).$$\\
(2) For any $\alpha\geq 1$,
\begin{equation*}
\begin{split}
\|\partial_t^{\alpha+1}\Psi(t)\|_{W_0} \leq C_{\mu_0,\alpha}&
e^{\sqrt{\mu_0}
t}(\|\Psi_t(0)\|_{W_0}+\|\Psi(0)\|_{W_0}+\sqrt{P_{\Psi}(0)})\\
+
C&_{\mu_0,\alpha}\sum_{i=1}^{\alpha}(\|\partial_t^{i+1}\Psi(0)\|_{W_0}+
\sqrt{P_{\partial_t^{i}\Psi}(0)}).
\end{split}
\end{equation*}
\end{lem}

Proof. Take the inner product of (4.3) with $\Psi_t$. Then we get

$$(W_0 \Psi_{tt}, \Psi_t)=(L_0 \Psi, \Psi_t)$$
$$\Longleftrightarrow\frac{d}{dt}(W_0 \Psi_{t}, \Psi_t)=\frac{d}{dt}(L_0 \Psi,
\Psi).$$\

The above equivalence comes from the self-adjointness of
$L_0$. Next integrate the above with respect to $t$ to get

\begin{equation}
(W_0
\Psi_t(t),\Psi_t(t))=(L_0\Psi(t),\Psi(t))+(W_0\Psi_t(0),\Psi_t(0))
-(L_0\Psi(0), \Psi(0)).
\end{equation}\

Since  $(L_0 \Psi(t),\Psi(t)) \leq \mu_0 (W_0 \Psi(t), \Psi(t))$
for all $t$ and $-(L_0\Psi(0), \Psi(0))\leq P_{\Psi}(0)$, from
(4.4) we get

\begin{equation}
\|\Psi_t(t)\|_{W_0}^2 \leq \mu_0 \|\Psi(t)\|_{W_0}^2 +
\|\Psi_t(0)\|_{W_0}^2 + P_{\Psi}(0).
\end{equation}\

Since $\|\Psi(t)\|_W \leq \int_0^t \|\Psi_t(\tau)\|_W d\tau +
\|\Psi(0)\|_W$, plugging this into (4.5), we get

\begin{equation*}
\|\Psi_t(t)\|_{W_0} \leq \sqrt{\mu_0} \int_0^t
\|\Psi_t(\tau)\|_{W_0}d\tau +
C(\|\Psi_t(0)\|_{W_0}+\|\Psi(0)\|_{W_0}+\sqrt{P_{\Psi}(0)}).
\end{equation*}\

By Gronwall's inequality, we obtain
\begin{equation*}
\begin{split}
\|\Psi_t(t)\|_{W_0} &\leq C e^{\sqrt{\mu_0}
t}(\|\Psi_t(0)\|_{W_0}+\|\Psi(0)\|_{W_0}+\sqrt{P_{\Psi}(0)})
\text{ and in success }\\
\
\|\Psi(t)\|_{W_0} &\leq C e^{\sqrt{\mu_0}  t}(
\|\Psi_t(0)\|_{W_0}+\|\Psi(0)\|_{W_0}+\sqrt{P_{\Psi}(0)}).
\end{split}
\end{equation*}\

Notice that $C$ only depends on $\mu_0$. For higher derivatives,
take $\partial_t^{\alpha}$ of (4.3):
$W_0(r)\partial_t^{\alpha}\Psi_{tt} = L_0\partial_t^{\alpha}
\Psi$. Take the inner product of this with $\partial_t
\partial_t^{\alpha} \Psi$ to get

\begin{equation*}
\begin{split}
\|\partial&_t^{\alpha+1}\Psi(t)\|_{W_0}^2 \\&=(W_0
\partial_t^{\alpha+1}\Psi(t),\partial_t^{\alpha+1}\Psi(t))\\
&=(L_0\partial_t^{\alpha}\Psi(t),\partial_t^{\alpha}\Psi(t))+
(W_0\partial_t^{\alpha+1}\Psi_t(0),\partial_t^{\alpha+1}\Psi_t(0))
-(L_0\partial_t^{\alpha}\Psi(0), \partial_t^{\alpha}\Psi(0))\\
&\leq \mu_0\|\partial_t^{\alpha}\Psi(t)\|_{W_0}^2 +
\|\partial_t^{\alpha+1}\Psi(0)\|_{W_0}^2 +
P_{\partial_t^{\alpha}\Psi}(0).
\end{split}
\end{equation*}\

Thus (2) easily follows.$\square$\\

Next we show that the energy estimates with Lemma 4.1 lead to the
same exponential growth rate on the $\sigma$ satisfying (4.1) and
(4.2). To avoid the confusion with the nonlinear analysis, we use
$\Phi$ instead of $\sigma$.

Define a weight function for $\Phi$ by
$$V_0(r)\equiv r^2\rho_0^{-\frac{4}{5}}=r^2(1+r^2)^2.$$

Notice that $(V_0, W_0)$ is chosen not randomly but to be a
symmetrizer of (4.1) and (4.2) that makes the energy estimate
work. $\|\Phi\|_{V_0}+\|\Psi\|_{W_0}$ resembles the real energy
(1.7).\\

\begin{lem} There are constants $C_1, C_2>0$ such
that\
$$\|\partial_t^{\alpha}\Phi(t)\|_{V_0}\leq C_1
\|\partial_t^{\alpha-1}\Psi(t)\|_{W_0}+\|\partial_t^{\alpha}
\Phi(0)\|_{V_0} +\|\partial_t^{\alpha}\Psi(0)\|_{W_0}, \text{ for
any }\alpha\geq 1$$\
 and $\|\Phi(t)\|_{V_0}\leq C_2 e^{\sqrt{\mu_0}
t}(\|\Phi_t(0)\|_{V_0}+\|\Phi(0)\|_{V_0}+\|\Psi_t(0)\|_{W_0}+
\|\Psi(0)\|_{W_0}+\sqrt{P_{\Psi}(0)}).$
\end{lem}

Proof. Fix $\alpha\geq 1$. Compute the following equation:

$$\int_0^{\infty}[\rho_0^{-\frac{4}{5}}\partial_t^{\alpha}\Phi\cdot
\partial_t^{\alpha}(4.1)+\frac{15}{4\pi}\rho_0\partial_t^{\alpha}\Psi
\cdot\partial_t^{\alpha}(4.2)] r^2 dr=0.$$

The choice of weight functions $\rho_0^{-\frac{4}{5}}$ and
$\frac{15}{4\pi}\rho_0$ above yields a nice cancellation after
integrating by parts, i.e.
$\int_0^{\infty}\rho_0^{-\frac{4}{5}}\partial_t^{\alpha}\Phi\cdot(r^2\rho_0
\partial_t^{\alpha}\Psi)_r+r^2\rho_0\partial_t^{\alpha}\Psi\cdot(\rho_0^
{-\frac{4}{5}}\partial_t^{\alpha}\Phi)_r dr=0$ and it results in

\begin{equation}
\frac{1}{2}\frac{d}{dt}\int_0^{\infty}[\rho_0^{-\frac{4}{5}}
(\partial_t^{\alpha}\Phi)^2+\frac{15}{4\pi}\rho_0
(\partial_t^{\alpha}\Psi)^2]r^2dr+15\int_0^{\infty}\rho_0\partial_t^
{\alpha}\Psi(\int_0^r\partial_t^{\alpha}\Phi s^2 ds)dr=0.
\end{equation}

On the other hand, by (4.1), we get
$\int_0^r\partial_t^{\alpha}\Phi s^2
ds=-r^2\rho_0\partial_t^{\alpha-1}\Psi$ and hence (4.6) can be
rewritten as
\begin{equation}
\frac{1}{2}\frac{d}{dt}[\int_0^{\infty}\rho_0^{-\frac{4}{5}}
(\partial_t^{\alpha}\Phi)^2r^2
dr+\int_0^{\infty}\frac{15}{4\pi}\rho_0
(\partial_t^{\alpha}\Psi)^2r^2
dr]=\frac{15}{2}\frac{d}{dt}\int_0^{\infty}
\rho_0^2(\partial_t^{\alpha-1}\Psi)^2 r^2 dr.
\end{equation}

As taking $\int_0^t$ of (4.7), the desired result is obtained.
Note that we have used $\rho_0^2\leq \rho_0$. As for $\alpha=0$,
utilize $\|\Phi(t)\|_{V_0}\leq \int_0^t
\|\Phi_t(\tau)\|_{V_0}d\tau
+\|\Phi(0)\|_{V_0}$.$\square$\\

The next two lemmas show that $\mu_0$ also determines the
exponential growth rate even with strong weights. Lemma 4.3 will play a
crucial role in the proof of Theorem 1.1.  Main idea of proofs
is to utilize the linear operator $L$. The results only contain
$\Psi$ estimates and the estimates on $\Phi$ can be derived similarly.\\

\begin{lem} For any $\alpha\geq 0$, there exists $C_{\mu_0}>0$
such that\
$$\int_0^{\infty}
\frac{r^2}{(1+r^2)^3}(\partial_t^{\alpha}\Psi_r)^2 dr + 2
\int_0^{\infty}\frac{1}{(1+r^2)^4} (\partial_t^{\alpha}\Psi)^2 dr
\leq C_{\mu_0} e^{2\sqrt{\mu_0} t}I_{\alpha},$$\ where
$I_{\alpha}=\sum_{i=0}^{\alpha+2}\|\partial_t^i\Psi(0)\|_{W_0}^2+
\sum_{i=0}^{\alpha+1}P_{\partial_t^{i}\Psi}(0)$ is given initial
data, for every solution $\Psi$ to (4.3).
\end{lem}

Proof. Multiply $\partial_t^{\alpha}$(4.3) by
$\partial_t^{\alpha}\Psi$ and integrate to get
\begin{equation}
\int_0^{\infty}W_0\partial_t^{\alpha+2}\Psi\partial_t^{\alpha}\Psi
dr = -\int_0^{\infty}
\frac{r^2}{(1+r^2)^3}(\partial_t^{\alpha}\Psi_r)^2 dr +
2\int_0^{\infty}\frac{3r^2-1}
{(1+r^2)^4}(\partial_t^{\alpha}\Psi)^2 dr.
\end{equation}\

If we move two negative terms in the RHS of (4.8) into the LHS, we
obtain
\begin{equation*}
\begin{split}
\int_0^{\infty}&
\frac{r^2}{(1+r^2)^3}(\partial_t^{\alpha}\Psi_r)^2 dr + 2
\int_0^{\infty}\frac{1}{(1+r^2)^4} (\partial_t^{\alpha}\Psi)^2
dr \\
&=6\int_0^{\infty}\frac{r^2}
{(1+r^2)^4}(\partial_t^{\alpha}\Psi)^2 dr
-\int_0^{\infty}W_0\partial_t^{\alpha+2}\Psi\partial_t^{\alpha}\Psi
dr \\
&\leq 6\int_0^{\infty}\frac{r^2}
{(1+r^2)^4}(\partial_t^{\alpha}\Psi)^2 dr +
\frac{1}{2}(\|\partial_t^{\alpha+2}\Psi(t)\|_{W_0}^2 +
\|\partial_t^{\alpha}\Psi(t)\|_{W_0}^2)\\
&\leq C (\|\partial_t^{\alpha+2}\Psi(t)\|_{W_0}^2 +
\|\partial_t^{\alpha}\Psi(t)\|_{W_0}^2)\\
\end{split}
\end{equation*}

But the last quantity is bounded by $C e^{2\sqrt{\mu_0}
t}I_{\alpha}$ by Lemma 4.1 and the conclusion follows.$\square$\\

For the next lemma define the weight functions $W_l$ for $l \leq
0$ by $$W_l(r)\equiv (1+r^2)^{-\frac{5+5l}{2}}W_0(r)=
\frac{15}{4\pi}r^2\rho_0^{1+l}.$$

We remark that $W_l$ is a linear version of symmetrizers $S_l$
defined in Section 5.\\

\begin{lem}For any $\alpha\geq0$, there exists $C_{\mu_0}>0$ such
that\
 $$\|\partial_t^{\alpha}\Psi(t)\|_{W_k} \leq C_{\mu_0}
e^{\sqrt{\mu_0} t}\sqrt{I_{\alpha}},$$ where $I_{\alpha}$ is
defined in Lemma 4.3, for every solution $\Psi$ to (4.3).
\end{lem}

Proof. Let us only consider $\alpha=0,1$. Other cases can be
treated in the similar way. Letting
$L_k=(1+r^2)^{\frac{k}{2}}L_0$, we get $(W_k\Psi_{tt},
\Psi_t)=(L_k \Psi, \Psi_t)$. But because $L_k$ is not self-adjoint
any more, we do not have a simple equivalent expression as in
Lemma 4.1. Instead, we try to find a relation between $L_k$ and
$L_0$. The following identities are needed and they are obtained
by definitions and straightforward computations.
$$(\Psi(1+r^2)^{\frac{k}{4}})_r = (1+r^2)^{\frac{k}{4}}\Psi_r +
\frac{k}{2}r(1+r^2)^{\frac{k}{4}-1}\Psi$$
$$(L_k \Psi, \Psi)=\int_0^{\infty}(1+r^2)^{\frac{k}{2}}\Psi
(\frac{r^2}{(1+r^2)^3}\Psi_r)_r dr +
2\int_0^{\infty}(1+r^2)^{\frac{k}{2}}\frac{3r^2-1}{(1+r^2)^4}\Psi^2
dr $$
\begin{equation*}
\begin{split}
(L_0(\Psi(1+r^2)^{\frac{k}{4}}), &\Psi(1+r^2)^{\frac{k}{4}})= -
\int_0^{\infty}\frac{r^2}{(1+r^2)^3}\{(1+r^2)^{\frac{k}{2}}\Psi_r^2
+k r (1+r^2)^{\frac{k}{2}-1}\Psi\Psi_r  \\
& +\frac{k^2}{4} r^2 (1+r^2)^{\frac{k}{2}-2}\Psi^2\}dr +2
\int_0^{\infty}\frac{3r^2-1}{(1+r^2)^4}(1+r^2)^{\frac{k}{2}}\Psi^2
dr \\
\
\end{split}
\end{equation*}
The last two identities imply the following:
$$(L_k \Psi, \Psi)=
(L_0(\Psi(1+r^2)^{\frac{k}{4}}), \Psi(1+r^2)^{\frac{k}{4}}) +
\frac{k^2}{4}\int_0^{\infty}\frac{r^4}{(1+r^2)^5}(\Psi
(1+r^2)^{\frac{k}{4}})^2 dr.$$

Define $Q_k$ and $I_k$ similar to $Q$ and $I$ in Section 2 as
following
\begin{equation*}
\begin{split}
Q_k(\Psi)&\equiv (L_0(\Psi(1+r^2)^{\frac{k}{4}}),
\Psi(1+r^2)^{\frac{k}{4}})-\beta
\int_0^{\infty}\frac{r^4}{(1+r^2)^5}(\Psi
(1+r^2)^{\frac{k}{4}})^2 dr ,\\
I_k(\Psi)&\equiv I(\Psi(1+r^2)^{\frac{k}{4}})=(W_k\Psi, \Psi).\\
\end{split}
\end{equation*}
where $\beta$ is a small positive constant. Then by doing the same
variational analysis in Section 2, one can show that there exists $\mu_k>0$ such
that $\mu_k$ is the maximum of a functional
$\frac{Q_k(\Psi)}{I_k(\Psi)}$ and hence $Q_k(\Psi)\leq\mu_k
I_k(\Psi)$.\\

Claim. $\mu_k < \mu_0$.\\

To see the Claim, pick a $\Psi_1$ such that $Q_k(\Psi_1)=\mu_k$
and $I_k(\Psi_1)=1$. By the definition of $Q_k$,
$$Q_k(\Psi_1)=
(L_0(\Psi_1(1+r^2)^{\frac{k}{4}}),
\Psi_1(1+r^2)^{\frac{k}{4}})-\beta
\int_0^{\infty}\frac{r^4}{(1+r^2)^5}(\Psi_1
(1+r^2)^{\frac{k}{4}})^2 dr.$$ Since
$(L_0(\Psi_1(1+r^2)^{\frac{k}{4}}),
\Psi_1(1+r^2)^{\frac{k}{4}})\leq \mu_0
I(\Psi_1(1+r^2)^{\frac{k}{4}})=\mu_0 I_k(\Psi_1)=\mu_0$,
$$Q_k(\Psi_1)=\mu_k \leq \mu_0 - \beta
\int_0^{\infty}\frac{r^4}{(1+r^2)^5}(\Psi_1
(1+r^2)^{\frac{k}{4}})^2 dr < \mu_0.$$
Thus,
\begin{equation}
\begin{split}
(L_k\Psi,\Psi)&=Q_k(\Psi) +
(\beta+\frac{k^2}{4})\int_0^{\infty}\frac{r^4}{(1+r^2)^5}(\Psi
(1+r^2)^{\frac{k}{4}})^2 dr \\
&\leq \mu_k (W_k\Psi,\Psi)
+(\beta+\frac{k^2}{4})\int_0^{\infty}\frac{r^4}{(1+r^2)^{\frac{10-k}{2}}}
\Psi^2 dr \\
\end{split}
\end{equation}

Now we are ready to go back to $(W_k\Psi_{tt}, \Psi_t)=(L_k \Psi,
\Psi_t)$.
\begin{equation*}
\begin{split}
&(L_k \Psi, \Psi_t)=\int_0^{\infty}(\frac{r^2}{(1+r^2)^3}\Psi_r)_r
(1+r^2)^{\frac{k}{2}}\Psi_t dr
+2\int_0^{\infty}\frac{3r^2-1}{(1+r^2)^{4}}(1+r^2)^{\frac{k}{2}}
\Psi\Psi_t dr \\
&=-\int_0^{\infty}\frac{r^2}{(1+r^2)^{\frac{6-k}{2}}}\Psi_r\Psi_{rt}dr
-k\int_0^{\infty}\frac{r^3}{(1+r^2)^{\frac{8-k}{2}}}\Psi_r\Psi_t
dr +
2\int_0^{\infty}\frac{3r^2-1}{(1+r^2)^{\frac{8-k}{2}}}\Psi\Psi_t
dr \\
&=\frac{1}{2}\frac{d}{dt}(L_k\Psi,\Psi)-k
\int_0^{\infty}(\frac{r^3}{(1+r^2)^{\frac{8-k}{2}}})_r\Psi\Psi_t
dr
-k\int_0^{\infty}\frac{r^3}{(1+r^2)^{\frac{8-k}{2}}}\Psi_r\Psi_t
dr \\
\end{split}
\end{equation*}

Hence,
\begin{equation}
\begin{split}
\frac{1}{2}\frac{d}{dt}(W_k\Psi_t,
\Psi_t)=\frac{1}{2}\frac{d}{dt}&(L_k\Psi,\Psi)-k
\int_0^{\infty}(\frac{r^3}{(1+r^2)^{\frac{8-k}{2}}})_r\Psi\Psi_t
dr\\
&-k\int_0^{\infty}\frac{r^3}{(1+r^2)^{\frac{8-k}{2}}}\Psi_r\Psi_t
dr.
\end{split}
\end{equation}

Let $k=1$. By Lemma 4.1 and 4.3, using the Cauchy-Schwartz
inequality,
$$-\int_0^{\infty}(\frac{r^3}{(1+r^2)^{\frac{7}{2}}})_r\Psi\Psi_t
dr -\int_0^{\infty}\frac{r^3}{(1+r^2)^{\frac{7}{2}}}\Psi_r\Psi_t
dr \leq C e^{2\sqrt{\mu_0} t} I_0$$ for a constant $C$ and
initial data $I_0$. Rewriting the above (4.10) when $k=1$, get
$$\frac{1}{2}\frac{d}{dt}(W_1\Psi_t,
\Psi_t)\leq \frac{1}{2}\frac{d}{dt}(L_1\Psi,\Psi)+ C
e^{2\sqrt{\mu_0} t} I_0.$$

Taking the integral with respect to $t$ and using (4.9), we have
\begin{equation*}
\begin{split}
(W_1\Psi_t,\Psi_t)&\leq (L_1\Psi,\Psi) +C e^{2\sqrt{\mu_0} t} I_0\\
&\leq \mu_1(W_1\Psi,
\Psi)+(\beta+\frac{1}{4})\int_0^{\infty}\frac{r^4}{(1+r^2)^{\frac{9}{2}}}
\Psi^2 dr +C e^{2\sqrt{\mu_0} t} I_0\\
\end{split}
\end{equation*}

Using Lemma 4.1 again, we obtain
\begin{equation}
(W_1\Psi_t,\Psi_t)\leq\mu_1(W_1\Psi,\Psi)+Ce^{2\sqrt{\mu_0} t}I_0.
\end{equation}

Since $\|\Psi(t)\|_W \leq \int_0^t \|\Psi_t(\tau)\|_W d\tau +
\|\Psi(0)\|_W$, combining this with (4.11), get
$$\|\Psi(t)\|_{W_1} \leq \sqrt{\mu_1}\int_0^{t}
\|\Psi(\tau)\|_{W_1} d\tau + Ce^{\sqrt{\mu_0} t}\sqrt{I_0}.$$

Since we have $\sqrt{\mu_1}<\sqrt{\mu_0}$, Gronwall inequality
gives $\|\Psi(t)\|_{W_1}\leq C e^{\sqrt{\mu_0}}\sqrt{I_0}$ as
well as $\|\Psi_t(t)\|_{W_1}\leq C e^{\sqrt{\mu_0} t}\sqrt{I_0}$.
The standard induction on $k$ with (4.10) claims the desired
result for all $k$.$\square$\\

\section{Weighted Instant Energy and Weighted Total Energy}

In this section, by the utilization of symmetrizers of the
Euler-Poisoon system, we introduce suitable measurements of
perturbations $\sigma$ and $u$ of steady states resembling
weighted Sobolev norms: the weighted instant energy
$\mathcal{E}_l$ and the weighted total energy
$\widetilde{\mathcal{E}}_l$, where $l$ is an index associated to
weights. The symmetrizers will play the same role as weights in
the linear analysis for small solutions. Then the total energy is
shown to be bounded by the instant energy under a certain
smallness assumption in using the equations directly, which makes
it sufficient to play only with the instant energy. The weighted
energy estimates for the instant energy will be carried out in the
next section. Before going any farther we remark that it is
convenient to work on rectangular coordinates rather than polar
coordinates because one can avoid the singularity of the origin
coming from the spherical symmetry.

We are interested in sufficiently small solutions $\sigma$, $u$
 satisfying the neutrality condition $$\int_{\mathbb{R}^3}\sigma dx=0.$$
$\sigma$ is assumed to be relatively smaller than $\rho_0$, in
particular, we assume
\begin{equation}
\frac{9}{10}\rho_0 \leq \rho_o+\sigma \leq \frac{11}{10}\rho_0.
\end{equation}
For such small solutions the Euler-Poisson system (1.1), (1.2) and
(1.3) can be rewritten in the rectangular coordinates as the
following:

\begin{equation}
\sigma_t+(\rho_0+\sigma)\nabla\cdot u +
\nabla(\rho_0+\sigma)\cdot u=0
\end{equation}
\begin{equation}
u_t + (u\cdot\nabla)u +\frac{4\pi}{15}(\rho_0+\sigma)
^{-\frac{4}{5}}\nabla\sigma -\{\frac{16\pi}{75}
\rho_0^{-\frac{9}{5}}\sigma+h(\sigma,\rho_0)
\}\nabla\rho_0+\nabla\phi=0
\end{equation}
\begin{equation}
\triangle\phi=4\pi\sigma
\end{equation}
where $h(\sigma,\rho_0)=-\frac{4\pi}{15}\{(\rho_0+\sigma)^
{-\frac{4}{5}}-\rho_0^{-\frac{4}{5}}+\frac{4}{5}\rho_0^
{-\frac{9}{5}}\sigma\}$ represents higher order terms. $u$ takes
the vector form, i.e. $u(x,t)=u(r,t)\frac{x}{r}$, where
$x\in\mathbb{R}^3$, $r=|x|$ and we denote each component of $u$
by $u^k$. Note $\nabla\times u=0$.

Now let us consider the symmetrizers $S_l$ where $l\in\mathbb{R}$
for the Euler-Poisson system.

\[ S_l = \begin{pmatrix}
         \frac{4\pi}{15}(\rho_0+\sigma)^{-\frac{4}{5}+l} & 0 & 0 &0\\
         0 & (\rho_0+\sigma)^{1+l} & 0 & 0 \\
         0 & 0 &(\rho_0+\sigma)^{1+l} & 0\\
         0 & 0 & 0 &(\rho_0+\sigma)^{1+l}\\
         \end{pmatrix}
\]\\

Define the instant energy $\mathcal{E}_l(t)$ and the total energy
$\widetilde{\mathcal{E}_l}(t)$ by
\begin{equation*}
\begin{split}
\mathcal{E}_l(t)\equiv&\sum_{j=0}^{3}\int_{\mathbb{R}^3}S_l
(\partial_t^j\sigma, \partial_t^j u^1, \partial_t^j u^2,
\partial_t^j u^3)^t\cdot(\partial_t^j\sigma, \partial_t^j u^1, \partial_t^j u^2,
\partial_t^j u^3)dx\\
=&\sum_{j=0}^{3}\int_{\mathbb{R}^3}\frac{4\pi}{15}
(\rho_0+\sigma)^{-\frac{4}{5}+l}(\partial_t^j\sigma)^2  +
(\rho_0+\sigma)^{1+l}|\partial_t^j u|^2dx\\
\equiv&\sum_{j=0}^{3}\mathcal{E}_l^j\\
\
\widetilde{\mathcal{E}_l}(t)\equiv&\sum_{j=0}^{3}\sum_{i=0}^j
\int_{\mathbb{R}^3}\frac{4\pi}{15}
(\rho_0+\sigma)^{-\frac{4}{5}+l+\frac{i}{5}}|\partial_t^{j-i}
\partial_x^{i}\sigma|^2 + (\rho_0+\sigma)^
{1+l+\frac{i}{5}}
|\partial_t^{j-i}\partial_x^{i} u|^2 dx\\
\equiv&\sum_{j=0}^{3}\sum_{i=0}^j\widetilde{\mathcal{E}_l}^{j,i}.
\end{split}
\end{equation*}\

Here $\partial_x$ represents any spatial first derivatives. Note
that $\mathcal{E}_0^0(t)$ is a part of the real energy 2$E$ of
$(\rho, u)$ defined in (1.7). The case $l=0$, however, is not
enough for our purpose because we cannot close the energy
estimate at the step $l=0$. Being convinced by Lemma 4.4,
we try different $l$'s as well. As
$l<0$ gets smaller, the weights become stronger due to the behavior
of $\rho_0$. Unfortunately,
as $l$ varies, new quadratic terms come out while performing the
energy estimates. This phenomenon seems undesirable but it turns
out that they are equipped with weaker weights. And that opens
another door.

Observe that the weights of mixed derivative terms in
$\widetilde{\mathcal{E}_l}$ are different, in fact a little better,
from the ones of temporal derivative terms. This is not a coincidence
but rather a nature of the system; the same feature can be seen
in the linear analysis. See Lemma 4.3. $\widetilde{\mathcal{E}_l}$
contains all the spatial and mixed derivatives of $\sigma$, $u$
and it is easy to see that
$$\mathcal{E}_l=\sum_{j=0}^3\widetilde{\mathcal{E}_l}^{j,0}
\leq\widetilde{\mathcal{E}_l}.$$ Now we want to show the converse,
in other words, $\widetilde{\mathcal{E}_l}$ is also bounded by
$\mathcal{E}_l$ under a certain smallness assumption:
\begin{equation}
\widetilde{\mathcal{E}_{-\frac{6}{5}}}+
\widetilde{\mathcal{E}_{-\frac{7}{5}}}+
\widetilde{\mathcal{E}_{-\frac{8}{5}}}\leq\theta_1
\end{equation}
where $\theta_1$ is a sufficiently small constant. In order to appreciate
the utilization of such an assumption, first we prove the next
lemma.\\

Notation. $\int dx$ represents $\int_{\mathbb{R}^3}dx$ and when
dealing with line integrals $\int dr$, each end value will be
specified.

\begin{lem} Suppose (5.5) holds for $0\leq t \leq T$. Then there
exists a constant $C>0$ such that for each $0\leq t \leq T$,
\begin{equation}
\sup_{x\in\mathbb{R}^3}|\frac{\sigma}{\rho_0+\sigma}|+
|\frac{\sigma_t}{\rho_0+\sigma}|+
|\frac{\nabla\sigma}{(\rho_0+\sigma) ^{\frac{9}{10}}}|+
|\frac{u}{(\rho_0+\sigma)^{\frac{1}{10}}}|+
|\frac{u_t}{(\rho_0+\sigma)^{\frac{1}{10}}}|+|\nabla u|\leq
C\sqrt{\theta_1}.
\end{equation}
In particular, the assumption (5.1) is justified.
\end{lem}

Proof. The above smallness assumption (5.5) together with the
Sobolev imbedding theorem yields the result. To see how it works,
let us apply the Sobolev imbedding theorem to
$|\frac{\sigma}{\rho_0+\sigma}|$.
\begin{equation}
\begin{split}
\sup|\frac{\sigma}{\rho_0+\sigma}|^2\leq&C
\|\frac{\sigma}{\rho_0+\sigma}\|_{H^2({\mathbb{R}^3})}^2\\
\leq&C\{\int|\frac{\sigma}{\rho_0+\sigma}|^2dx+
\int|\frac{\partial_x\sigma}{\rho_0+\sigma}|^2dx+
\int|\frac{\partial_x^2\sigma}{\rho_0+\sigma}|^2dx\\
&+\int|\frac{\sigma}{\rho_0+\sigma}\frac{\partial_x(\rho_0+\sigma)}
{\rho_0+\sigma}|^2dx+ \int|\frac{\partial_x\sigma}{\rho_0+\sigma}
\frac{\partial_x(\rho_0+\sigma)}{\rho_0+\sigma}|^2dx\\
&+\int|\frac{\sigma}{\rho_0+\sigma}\frac{\partial_x^2(\rho_0+\sigma)}
{\rho_0+\sigma}|^2dx
+\int\{|\frac{\sigma}{\rho_0+\sigma}||\frac{\partial_x(\rho_0+\sigma)}
{\rho_0+\sigma}|^2\}^2dx\}
\end{split}
\end{equation}
Recall the behavior of $\rho_0$, namely $\rho_0(r)=O(r^{-5})$ for
large $r$. Hence $|\frac{\partial_x\rho_0}{\rho_0+\sigma}|$ and
$|\frac{\partial_x^2\rho_0}{\rho_0+\sigma}|$ are uniformly
bounded. Thus (5.7) becomes
\begin{equation}
\begin{split}
\sup|\frac{\sigma}{\rho_0+\sigma}|^2\leq &C\{(\text{
}\widetilde{\mathcal{E} _{-\frac{6}{5}}}^{0,0}+\widetilde{\mathcal{E}
_{-\frac{7}{5}}}^{1,1}+\widetilde{\mathcal{E}
_{-\frac{8}{5}}}^{2,2}\text{ })\\
&+\sup|\frac{\sigma}{\rho_0+\sigma}|^2 (\text{
}\widetilde{\mathcal{E} _{-\frac{6}{5}}}^{0,0}+\widetilde{\mathcal{E}
_{-\frac{7}{5}}}^{1,1}+\widetilde{\mathcal{E}
_{-\frac{8}{5}}}^{2,2}\text{ })\\
&+\sup|\frac{\sigma}{\rho_0+\sigma}|^2
\sup|\frac{\partial_x\sigma}{(\rho_0+\sigma)^{\frac{9}{10}}}|^2
\text{ }\widetilde{\mathcal{E}_{-\frac{7}{5}}}^{1,1}\text{ }\}\\
\leq&C(\theta_1+\sup|\frac{\sigma}{\rho_0+\sigma}|^2\theta_1+
\sup|\frac{\sigma}{\rho_0+\sigma}|^2
\sup|\frac{\partial_x\sigma}{(\rho_0+\sigma)^{\frac{9}{10}}}|^2
\theta_1).
\end{split}
\end{equation}

We can get similar estimates to (5.8) for other terms in (5.6).
Call the LHS of (5.6) $S$. Consequently, we obtain the following:
\begin{equation}
S^2\leq C\theta_1+CS^2\theta_1+CS^4\theta_1
\end{equation}
Since $\theta_1$ is small enough, (5.9) immediately implies the
lemma.$\square$\\

Since $\widetilde{\mathcal{E}_l}^{j,0}=\mathcal{E}_l^j$, we only
need to show that $\widetilde{\mathcal{E}_l}^{j,i}$ where $1\leq
i\leq j\leq 3$ is bounded by $C\mathcal{E}_l$; the precise
statement is given in the following lemma.  Main idea for
accomplishing the goal is to estimate the spatial and mixed
derivative terms directly from the equations in terms of temporal
derivative terms.\\

\begin{lem} Suppose (5.5) holds for $0\leq t \leq T$. Let
$l\leq 0$. Then there exists $C>0$ such that for each $1\leq i
\leq j \leq 3$ and for $0\leq t \leq T$,
$$\widetilde{\mathcal{E}_l}^{j,i}(t)\leq
C\sum_{k=1}^j\mathcal{E}_l^k(t)+C\sum_{k=0}^{j-1}\mathcal{E}_
{l+\frac{3}{5}}^k.$$
\end{lem}

Proof. Let us start with one spatial derivative terms
corresponding to $i=1$. No temporal derivative terms i.e. when
$i=1\text{, }j=1$ are treated carefully since it is instructive
and other cases can be easily shown from it. Various case of $j$
for a fixed $i$ will be done in turn. Then we move onto other cases
of $i$. First of all, solve (5.2) and (5.3) for $\nabla\sigma$ and
$\nabla\cdot u$ to get

\begin{equation}
\nabla\sigma=-\frac{15}{4\pi}(\rho_0+\sigma)^{\frac{4}{5}}\{
u_t-\frac{16\pi}{75}\rho_0^{-\frac{9}{5}}\nabla\rho_0\sigma
+\nabla\phi+(u\cdot\nabla)u- h \nabla\rho_0 \}
\end{equation}
\begin{equation}
\nabla\cdot u=-\frac{1}{\rho_0+\sigma}\{\sigma_t+\nabla\rho_0
\cdot u+\nabla\sigma\cdot u\}.
\end{equation}\

Notice that the estimate on $\nabla\cdot u$ is enough for
$\partial_x u$ since $\nabla\times u=0$. In order to get a right
weight of $\widetilde{\mathcal{E}_l}^{1,1}$ for $\sigma$ and $u$,
multiply (5.10) and (5.11) by
$(\rho_0+\sigma)^{-\frac{3}{10}+\frac{l}{2}}$ and
$(\rho_0+\sigma)^{\frac{3}{5}+\frac{l}{2}}$ respectively, square
and integrate them over $\mathbb{R}^3$.
\begin{equation}
\begin{split}
\int(\rho_0+&\sigma)^{-\frac{3}{5}+l}|\nabla\sigma|^2
dx +\int(\rho_0+\sigma)^{\frac{6}{5}+l}|\nabla\cdot u|^2dx\\
\leq C&\int
(\rho_0+\sigma)^{1+l}\{|u_t|^2+|\rho_0^{-\frac{9}{5}}\nabla\rho_0|^2
\sigma^2+|\nabla\phi|^2
+|(u\cdot\nabla)u|^2+|\nabla\rho_0|h^2\}dx\\
+&C\int(\rho_0+\sigma)^{-\frac{4}{5}+l}\{|\sigma_t|^2+|\nabla\rho_0
\cdot u|^2+|\nabla\sigma\cdot u|^2\}dx\\
\leq C&\int(\rho_0+\sigma)^{1+l}|u_t|^2+
(\rho_0+\sigma)^{-\frac{4}{5}+l}|\sigma_t|^2dx\\
+&C\int(\rho_0+\sigma)^{1+l}|\rho_0^{-\frac{9}{5}}\nabla\rho_0|^2
\sigma^2+(\rho_0+\sigma)^{-\frac{4}{5}+l}|\nabla\rho_0 \cdot
u|^2dx\\
+&C\int(\rho_0+\sigma)^{1+l}|(u\cdot\nabla)u|^2+
(\rho_0+\sigma)^{-\frac{4}{5}+l}|\nabla\sigma\cdot u|^2dx\\
+&C\int(\rho_0+\sigma)^{1+l}|\nabla\rho_0|h^2dx
+C\int(\rho_0+\sigma)^{1+l}|\nabla\phi|^2dx
\end{split}
\end{equation}
We rearranged terms according to the order of algebraic degree.
The key difficulty lies in the potential part and so it will be
done last. The first integral is bounded by $C\mathcal{E}_l^1$ by
the definition. Since
$|\rho_0^{-\frac{9}{5}}\nabla\rho_0|=5r(1+r^2)\leq C
(\rho_0+\sigma) ^{-\frac{3}{5}}$ and similarly $|\nabla\rho_0|\leq
C(\rho_0+\sigma)^{\frac{6}{5}}$,

\begin{equation}
\begin{split}
\int&(\rho_0+\sigma)^{1+l}|\rho_0^{-\frac{9}{5}}\nabla\rho_0|^2
\sigma^2 dx +\int(\rho_0+\sigma)^{-\frac{4}{5}+l}|\nabla\rho_0
\cdot u|^2 dx\\
&\leq C\int(\rho_0+\sigma)^{-\frac{1}{5}+l}\sigma^2dx+C\int
(\rho_0+\sigma)^{\frac{8}{5}+l}|u|^2 dx\\
&\leq C\mathcal{E}_{l+\frac{3}{5}}^0.
\end{split}
\end{equation}\

For higher order terms, we use the smallness assumption to get the
desired estimates: Lemma 5.1 is used.

\begin{equation}
\begin{split}
\int&(\rho_0+\sigma)^{1+l}|(u\cdot\nabla)u|^2+
(\rho_0+\sigma)^{-\frac{4}{5}+l}|\nabla\sigma\cdot u|^2dx\\
&\leq\int|\frac{u}{(\rho_0+\sigma)^{\frac{1}{10}}}|^2
\{(\rho_0+\sigma)^{\frac{6}{5}+l}|\nabla\cdot u|^2
+(\rho_0+\sigma)^{-\frac{3}{5}+l}|\nabla\sigma|^2\} dx\\
&\leq C\theta_1\int(\rho_0+\sigma)^{\frac{6}{5}+l}|\nabla\cdot
u|^2 +(\rho_0+\sigma)^{-\frac{3}{5}+l}|\nabla\sigma|^2 dx
\end{split}
\end{equation}\

Hence the higher order terms
in the third integral of the RHS in (5.12) can be absorbed into
its LHS, since $\theta_1$ is sufficiently small. Because $h$ is
a higher order term  depending on only
$\rho_0$ and $\sigma$, after applying Lemma 5.1 to the fourth
integral in (5.2), we get
\begin{equation}
\int(\rho_0+\sigma)^{1+l}|\nabla\rho_0|h^2dx\leq C\theta_1
\mathcal{E}_{l+\frac{3}{5}}^0.
\end{equation}\

In order to handle the potential part, we write the poisson
equation (5.4) for the spherically symmetric case in polar
coordinates as
$$\phi_{rr}+\frac{2}{r}\phi_r =4\pi\sigma \text{ where }
\phi_r=\frac{4\pi}{r^2}\int_0^r \sigma s^2 ds.$$ Recall the $L^p$
estimates $\|\partial^2\phi\|_{L^2(\mathbb{R}^3)} \leq
C\|\sigma\|_{L^2(\mathbb{R}^3)}$ and in the spherically symmetric
case it implies that since
$\sum_{i,j=1}^3(\partial_i\partial_j\phi)^2=\phi_{rr}^2
+\frac{2}{r^2}\phi_r^2$,
$$\|\frac{1}{r}\phi_r\|_{L^2(\mathbb{R}^3)}\leq C
\|\sigma\|_{L^2(\mathbb{R}^3)}.$$
We consider two cases:
$l\geq-\frac{3}{5}$ and $l<-\frac{3}{5}$. In the first case,
$(\rho_0+\sigma)^{1+l}r^2\sim(1+r^2)^{-\frac{5(1+l)}{2}}r^2$ is
uniformly bounded, and therefore we get
\begin{equation}
\int(\rho_0+\sigma)^{1+l}|\nabla\phi|^2dx=4\pi\int_0^{\infty}
(\rho_0+\sigma)^{1+l}\phi_r^2 r^2 dr\leq C \int_0^{\infty}
(\frac{\phi_r}{r})^2 r^2dr\leq C\int\sigma^2 dx.
\end{equation}
The $L^p$ estimate has been used at the last inequality. When
$l<-\frac{3}{5}$, we divide the integral into two parts.
\begin{equation}
\begin{split}
\int(\rho_0+\sigma)^{1+l}|\nabla\phi|^2dx&=4\pi\int_0^{\infty}
(\rho_0+\sigma)^{1+l}\phi_r^2 r^2
dr\\&=4\pi\int_0^{1}(\rho_0+\sigma)^{1+l} \phi_r^2
r^2dr+4\pi\int_1^{\infty}(\rho_0+\sigma)^{1+l}
\phi_r^2 r^2dr\\
&\equiv(I)+(II)\\
\end{split}
\end{equation}
In the unit ball we can do the same as in (5.16), since the weight
function is bounded.
\begin{equation}
(I)\leq C\int_0^1(\frac{\phi_r}{r})^2r^2 dr\leq C\int_0^1\sigma^2
r^2 dr\leq C \int (\rho_0+\sigma)^{m}\sigma^2 dx,\text{ for any
}m.
\end{equation}
For the second term $(II)$, we use the neutrality condition
$\int_{\mathbb{R}^3}\sigma dx=0$ which is equivalent to
$\int_0^r\sigma s^2 ds=-\int_r^{\infty}\sigma s^2 ds$ in polar
coordinates.
\begin{equation}
\begin{split}
(II)&\leq
C\int_1^{\infty}(1+r^2)^{-\frac{5(1+l)}{2}}(\frac{1}{r^2}\int_r^
{\infty}\sigma s^2 ds)^2 r^2 dr\\
&=C\int_1^{\infty}\frac{1}{r^2(1+r^2)}\{(1+r^2)^{-\frac{5}{4}
(\frac{3}{5}+l)}\int_r^{\infty}\sigma s^2 ds\}^2dr\\
&\leq C\int_1^{\infty}\frac{1}{r^2(1+r^2)}\{\int_r^{\infty}(1+s^2)
^{-\frac{5}{4} (\frac{3}{5}+l)}\sigma s^2 ds\}^2 dr\\
&\leq C\int_1^{\infty}\frac{1}{r^2(1+r^2)}[\int_r^{\infty}(1+s^2)
^{-\frac{5}{2}(-\frac{1}{5}+l)}\sigma^2 s^2 ds]
[\int_r^{\infty}\frac{s^2}{(1+s^2)^2}ds]dr\\
&\leq C[\int_1^{\infty}\frac{1}{r^2(1+r^2)}dr]
[\int_1^{\infty}\frac{s^2}{(1+s^2)^2}ds][\int_1^{\infty}(1+s^2)
^{-\frac{5}{2}(-\frac{1}{5}+l)}\sigma^2 s^2 ds]\\
&\leq C\int(\rho_0+\sigma)^{-\frac{1}{5}+l}\sigma^2 dx
\end{split}
\end{equation}

From (5.16), (5.18) and (5.19) we conclude that for any $l\leq 0$,
\begin{equation}
\int(\rho_0+\sigma)^{1+l}|\nabla\phi|^2dx\leq C
\mathcal{E}_{l+\frac{3}{5}}^0.
\end{equation}

Thus, from (5.12), (5.13), (5.14), (5.15) and (5.20), we obtain
$$\widetilde{\mathcal{E}_l}^{1,1}(t)\leq
C\mathcal{E}_l^1+C\mathcal{E}_{l+\frac{3}{5}}^0.$$\\

Next we focus on one spatial, one temporal derivative terms,
namely the case $i=1,\text{ }j=2$. Take $\partial_t$ of (5.10) and (5.11):
\begin{equation}
\begin{split}
\nabla\sigma_t&=-\frac{3}{\pi}\frac{\sigma_t}{\rho_0+\sigma}
(\rho_0+\sigma)^{\frac{4}{5}}\{
u_t-\frac{16\pi}{75}\rho_0^{-\frac{9}{5}}\nabla\rho_0\sigma
+\nabla\phi+(u\cdot\nabla)u- h \nabla\rho_0 \} \\
&-\frac{15}{4\pi}(\rho_0+\sigma)^{\frac{4}{5}}\{u_{tt}-\frac{16\pi}{75}
\rho_0^{-\frac{9}{5}}\nabla\rho_0\sigma_t
+\nabla\phi_t+(u_t\cdot\nabla)u+(u\cdot\nabla)u_t - h_t
\nabla\rho_0\}\\
\end{split}
\end{equation}
\begin{equation}
\begin{split}
\nabla\cdot
u_t&=\frac{\sigma_t}{\rho_0+\sigma}\frac{1}{\rho_0+\sigma}
\{\sigma_t+\nabla\rho_0 \cdot u +\nabla\sigma\cdot
u\}\\&-\frac{1}{\rho_0+\sigma}\{\sigma_{tt}+\nabla\rho_0 \cdot
u_t+\nabla\sigma_t\cdot u+\nabla\sigma\cdot u_t\}.
\end{split}
\end{equation}\

The first part of (5.21) and (5.22) is bounded by
$|\frac{\sigma_t}{\rho_0+\sigma}\nabla\sigma|$ and
$|\frac{\sigma_t}{\rho_0+\sigma}\nabla\cdot u|$ that have been
already estimated at the previous step. Note that
$|\frac{\sigma_t}{\rho_0+\sigma}|$ is small enough and it does
not cause any trouble. The other part has the same structure as
before and therefore we can do the same: multiply the same
weights used in the previous case, square and integrate. The
potential term is easily taken care of, since the dynamics of
$\nabla\phi_t$ gets simpler and better in the sense that
$$\nabla\phi_t=4\pi\nabla\triangle^{-1}\sigma_t=-
4\pi\nabla\triangle^{-1}\nabla\cdot(\rho_0+\sigma)u=-4\pi
(\rho_0+\sigma)u.$$ After higher order terms being absorbed we
obtain
$$\widetilde{\mathcal{E}_l}^{2,1}\leq C(\mathcal{E}_l^2+
\mathcal{E}_l^1)+C(\mathcal{E}_{l+\frac{3}{5}}^1+
\mathcal{E}_{l+\frac{3}{5}}^0).$$\

Considering $\partial_t(5.21)$ and $\partial_t(5.22)$, each term
in the RHS either has been estimated or can be dealt with in the
same manner as the previous cases, and therefore the estimates on
$\partial_t^2\partial_x\sigma$ and $\partial_t^2\partial_x u$
follows:
$$\widetilde{\mathcal{E}_l}^{3,1}\leq C(\mathcal{E}_l^3+\mathcal{E}_l^2+
\mathcal{E}_l^1)+C(\mathcal{E}_{l+\frac{3}{5}}^2+
\mathcal{E}_{l+\frac{3}{5}}^1+ \mathcal{E}_{l+\frac{3}{5}}^0).$$\\

Now we move onto two spatial derivative terms, the case $i=2$. The
only but important difference is to use a different weight to
close the estimate. It explains why the total energy is designed
with having different weights according to the number of spatial
derivatives. Compute $\partial_x (5.10)$ and $\partial_x (5.11)$
to get
\begin{equation}
\begin{split}
\nabla\partial_x\sigma=-\frac{3}{\pi}\frac{\partial_x(\rho_0+\sigma)}
{(\rho_0+\sigma)^{\frac{1}{5}}}\{&
u_t-\frac{16\pi}{75}\rho_0^{-\frac{9}{5}}\nabla\rho_0\sigma
+\nabla\phi+(u\cdot\nabla)u- h \nabla\rho_0 \}\\
-\frac{15}{4\pi}(\rho_0+\sigma)^{\frac{4}{5}}\{\partial_x
u_t&-\frac{16\pi}{75}
\rho_0^{-\frac{9}{5}}\nabla\rho_0\partial_x\sigma
+\nabla\partial_x\phi+(\partial_x
u\cdot\nabla)u+(u\cdot\nabla)\partial_x u\\
&- \partial_x h \nabla\rho_0 -\frac{16\pi}{75}
\partial_x(\rho_0^{-\frac{9}{5}}\nabla\rho_0)\sigma-h\nabla
\partial_x\rho_0\}\\
\end{split}
\end{equation}
\begin{equation}
\begin{split}
\nabla\cdot
\partial_x u&=\frac{\partial_x(\rho_0+\sigma)}
{(\rho_0+\sigma)^2} \{\sigma_t+\nabla\rho_0 \cdot u
+\nabla\sigma\cdot
u\}\\-\frac{1}{\rho_0+\sigma}\{&\partial_x\sigma_t+\nabla\rho_0
\cdot \partial_x u+\nabla\partial_x\sigma\cdot u+\nabla\sigma\cdot
\partial_x u +
\nabla\partial_x\rho_0 \cdot u\}.
\end{split}
\end{equation}

In order to get right exponents $-\frac{2}{5}+l$, $\frac{7}{5}+l $
of $\nabla\partial_x\sigma$ and $\nabla\cdot\partial_x u$ in
$\text{ }\widetilde{\mathcal{E}_l}^{2,2}$ we multiply (5.23) and
(5.24) by $(\rho_0+\sigma)^{-\frac{1}{5}+\frac{l}{2}}$ and
$(\rho_0+\sigma)^{\frac{7}{10}+\frac{l}{2}}$ respectively, and
square them. Notice that our chosen weight functions are of
polynomial type due to the behavior of $\rho_0$:
$(\rho_0+\sigma)\sim\rho_0=O(r^{-5})$. Thus as one takes the
spatial derivative, one gets $|\partial_x\rho_0|=O(r^{-6})\sim
(\rho_0+\sigma)^{\frac{6}{5}}$. So we get the following:
\begin{equation}
\begin{split}
\int(\rho_0+\sigma)^{-\frac{2}{5}+l}|\nabla\partial_x\sigma|^2
dx +\int(\rho_0+\sigma)^{\frac{7}{5}+l}|\nabla\cdot\partial_x& u|^2dx\\
\leq C\int ((\rho_0+\sigma)^{\frac{8}{5}+l}+
|\frac{\partial_x\sigma}{(\rho_0+\sigma)^{\frac{9}{10}}}|^2
(\rho_0+\sigma)^{1+l}
)(|u_t|&^2+|\rho_0^{-\frac{9}{5}}\nabla\rho_0|^2
\sigma^2+|\nabla\phi|^2\\
+|(u\cdot\nabla)u|^2+|\nabla\rho_0|^2h^2)dx\\
+C\int((\rho_0+\sigma)^{-\frac{1}{5}+l}+|\frac{\partial_x\sigma}
{(\rho_0+\sigma)^{\frac{9}{10}}}|^2
(\rho_0+\sigma)^{-\frac{4}{5}+l})(|&\sigma_t|^2+|\nabla\rho_0
\cdot u|^2+|\nabla\sigma\cdot u|^2)dx\\
+C\int (\rho_0+\sigma)^{\frac{6}{5}+l}(|\partial_x
u_t|^2+|\rho_0^{-\frac{9}{5}}\nabla\rho_0|^2
|\partial_x\sigma|^2+|\nabla\partial&_x\phi|^2 +|(\partial_x
u\cdot\nabla)u|^2\\
+|(u\cdot\nabla)\partial_x u|^2+|\nabla\rho_0|^2|\partial_x h|^2
+|\partial_x(\rho_0^{-\frac{9}{5}}\nabla\rho_0)&|^2\sigma^2
+|\nabla\partial_x\rho_0|^2h^2)dx\\
+C\int(\rho_0+\sigma)^{-\frac{3}{5}+l}(|\partial_x\sigma_t|^2
+|\nabla\rho_0 \partial_x u|^2+|\nabla\partial_x\sigma
u|^2&+|\nabla\sigma\partial_x\cdot u|^2+|\nabla\partial_x
\rho_0 u|^2)dx\\
\end{split}
\end{equation}

The first and second integrals are exactly same as the RHS of
(5.12). The ones in the third and fourth integrals except the
potential term have been already estimated at the previous steps
since each of them contains only one spatial derivative with the
right exponent of the corresponding weight. The potential part
does not produce any
further difficulty, indeed it behaves better both in weights and
in derivatives. Notice that
$\triangle\partial_x\phi=4\pi\partial_x\sigma$ and
$\int_{\mathbb{R}^3}\partial_x\sigma dx=0$. Thus we can do the
same as we did in (5.17), (5.18), (5.19) and we get the following
estimate similar to (5.20):
\begin{equation}
\int(\rho_0+\sigma)^{\frac{6}{5}+l}|\nabla\partial_x\phi|^2 dx
\leq C\mathcal{E}_{l+\frac{3}{5}}^1
\end{equation}

Consequently, we get the desired estimates for
$\partial_x^2\sigma$ and $\partial_x^2 u$:
\begin{equation*}
\widetilde{\mathcal{E}_l}^{2,2}\leq
C(\mathcal{E}_l^2+\mathcal{E}_l^1)+C(\mathcal{E}_{l+\frac{3}{5}}^1
+\mathcal{E}_{l+\frac{3}{5}}^0)
\end{equation*}\

Consider $\partial_t(5.23)$ and $\partial_t(5.24)$. Taking
$\partial_t$ does not destroy the structure of equations. As
going along the same track, the desired result on
$\partial_t\partial_x^2\sigma$ and $\partial_t\partial_x^2 u$ is
easily obtained.\\

Lastly, for three full spatial derivative terms, namely the case
$i=3\text{, }j=3$, compute $\partial_x(5.23)$ and
$\partial_x(5.24)$. Since we are dealing with one more spatial
derivative, we have to modify the weights again. Multiply them by
$(\rho_0+\sigma)^{-\frac{1}{10}+\frac{l}{2}}$ and
$(\rho_0+\sigma)^{\frac{4}{5}+\frac{l}{2}}$, square and integrate
them. Then most terms have been already treated before. As for
the potential part, noting that
$\triangle\partial_x^2\phi=4\pi\partial_x^2\sigma$ and
$\int_{\mathbb{R}^3}\partial_x^2\sigma dx=0$, we get the similar
estimate to (5.26):
$$\int(\rho_0+\sigma)^{\frac{7}{5}+l}|\nabla\partial_x^2\phi|^2 dx\leq
C\mathcal{E}_{l+\frac{3}{5}}^2$$ At last this finishes
the lemma.$\square$\\

Lemma 5.2 shows that any spatial and mixed derivative terms can
be estimated in terms of time derivative terms with suitable
weights, i.e. $\widetilde{\mathcal{E}_l}$ and $\mathcal{E}_l$ are
more or less equivalent measurements. Now we take time
derivatives which do not destroy the structure of the system much
and do the energy
estimates.\\

\section{Weighted Nonlinear Energy Estimates}

In this section we perform the nonlinear energy estimates with
the utilization of a family of symmetrizers of the system. Energy
estimates with weights on $\partial_t^j\sigma$ and $\partial_t^j
u$ for $0\leq j\leq 3$ are carried out to derive the following key
estimate so as to build the bootstrap argument which will be
discussed in the next section. Throughout this section, (5.1) and (5.5)
are assumed.\\

\begin{prop} Suppose $\widetilde{\mathcal{E}_{-\frac{6}{5}}}+
\widetilde{\mathcal{E}_{-\frac{7}{5}}}+
\widetilde{\mathcal{E}_{-\frac{8}{5}}}\leq\theta_1$
 for $0\leq t\leq T$ where $\theta_1\ll 1$ is sufficiently
 small.\\
(1) Let $l=0$. Then, for any fixed small $\eta>0$, there exist
$C$, $C_{\eta}>0$ such that\
\begin{equation*}
\text{ }
\frac{1}{2}\frac{d}{dt}\mathcal{E}_0\leq(C\sqrt{\theta_1}+\eta)\mathcal{E}_0
+C_{\eta}(\mathcal{E}_0)^{\frac{3}{2}}(\mathcal{E}_{-3})^{\frac{1}{2}}
 +C_{\eta}(\mathcal{E}_0^0+\mathcal{E}_0^1+\mathcal{E}_0^2).
\end{equation*}\\
(2) Let $l<0$. Then, for any fixed small $\eta>0$, there exist
$C$, $C_{\eta}>0$ such that\
\begin{equation*}
\text{ }\text{ }\text{ }\text{ }\text{
}\frac{1}{2}\frac{d}{dt}\mathcal{E}_l
\leq(C\sqrt{\theta_1}+\eta)\mathcal{E}_l+ C\mathcal{E}_{l+\frac{3}{10}}
+C_{\eta}(\mathcal{E}_l)^{\frac{3}{2}}(\mathcal{E}_{l+(-3-2l)})^{\frac{1}{2}}
+C_{\eta}(\mathcal{E}_l^0+\mathcal{E}_l^1+\mathcal{E}_l^2).
\end{equation*}\\
In particular, if $\text{ }l\leq-\frac{3}{2}$, then we have
\begin{equation*}
\text{ }\text{
}\frac{1}{2}\frac{d}{dt}\mathcal{E}_l\leq(C\sqrt{\theta_1}+\eta)\mathcal{E}_l+
C\mathcal{E}_{l+\frac{3}{10}} +C_{\eta}(\mathcal{E}_l)^2
+C_{\eta}(\mathcal{E}_l^0+\mathcal{E}_l^1+\mathcal{E}_l^2).
\end{equation*}
\end{prop}\

Proposition 6.1 will be proven by a series of lemmas in which we
will derive the estimate on $\mathcal{E}_l^j$ for each $j$ ; each lemma has its own significance and we will need
all of them to prove the bootstrap argument. Let us start with
the simplest case $j=0$: the zeroth order estimate.\\

\begin{lem}$(\mathcal{E}_l^0)$ Let $l\leq 0$. There exists a
constant $C>0$ such that
$$\frac{1}{2}\frac{d}{dt}\mathcal{E}_l^0 \leq C\sqrt{\theta_1}
\mathcal{E}_l^0 + C\mathcal{E}_{l+\frac{3}{10}}^0.$$
\end{lem}

 Proof. Consider $$\int S_l \binom{5.2}{5.3} \cdot
\binom{\sigma}{u} dx =0:$$

\begin{equation}
\begin{split}
0=&\frac{4\pi}{15}\int(\rho_0+\sigma)^{-\frac{4}{5}+l} \sigma_t
\sigma dx + \frac{4\pi}{15}\int(\rho_0+\sigma)^{-\frac{4}{5}+l}
[\nabla\cdot(\rho_0+\sigma)u]\sigma dx \\
&+\int(\rho_0+\sigma)^{1+l} u_t\cdot u dx +
\int(\rho_0+\sigma)^{1+l}(u \cdot\nabla)u\cdot u dx \\
&+ \frac{4\pi}{15}\int(\rho_0+\sigma)^{1+l}
[(\rho_0+\sigma)^{-\frac{4}{5}} \nabla\sigma -
\frac{4}{5}\rho_0^{-\frac{9}{5}}\nabla\rho_0\sigma +
h(\sigma, \rho_0)]\cdot u dx  \\
&+\int (\rho_0+\sigma)^{1+l}\nabla\phi\cdot u dx \\
\end{split}
\end{equation}

We compute the first three terms in turn; call them $(I)\text{,
}(II)$ and $(III)$.
\begin{equation}
\begin{split}
(I)=&\frac{1}{2}\frac{d}{dt}
\frac{4\pi}{15}\int(\rho_0+\sigma)^{-\frac{4}{5}+l} \sigma^2 dx
-\frac{1}{2}\frac{4\pi}{15}(-\frac{4}{5}+l)
\int(\rho_0+\sigma)^{-\frac{9}{5}+l}\sigma_t \sigma^2dx\\
(II)= &-\frac{4\pi}{15}(-\frac{4}{5}+l)\int(\rho_0+\sigma)^
{-\frac{9}{5}+l}\nabla(\rho_0+\sigma)\cdot [(\rho_0+\sigma)u]
\sigma dx\\&- \frac{4\pi}{15}\int(\rho_0+\sigma)^{-\frac{4}{5}+l}
(\rho_0+\sigma)u \cdot\nabla\sigma dx\\
(III)=&\frac{1}{2}\frac{d}{dt}\int(\rho_0+\sigma)^{1+l}|u|^2dx
-\frac {1}{2}(1+l)\int(\rho_0+\sigma)^{l}\sigma_t |u|^2dx
\end{split}
\end{equation}

After cancellation (6.1) becomes
\begin{equation}
\begin{split}
&\frac{1}{2}\frac{d}{dt}[\frac{4\pi}{15}\int(\rho_0+\sigma)^
{-\frac{4}{5}+l}\sigma^2 dx +
\int(\rho_0+\sigma)^{1+l}u^2 dx]\\
&=\frac{1}{2}\frac{4\pi}{15}(-\frac{4}{5}+l)
\int(\rho_0+\sigma)^{-\frac{4}{5}+l}\frac{\sigma_t}{\rho_0
+\sigma} \sigma^2 dx+\frac
{1}{2}(1+l)\int(\rho_0+\sigma)^{1+l}\frac{\sigma_t}
{\rho_0+\sigma} |u|^2 dx \\
&-\int(\rho_0+\sigma)^{1+l}(u \cdot\nabla)u\cdot u
dx-\int(\rho_0+\sigma)^{1+l}\nabla\phi\cdot u
dx \\
&+\frac{4\pi}{15}\int(\rho_0+\sigma)^{1+l}[(-\frac{4}{5}+l)
(\rho_0+\sigma)^{-\frac{9}{5}}\nabla(\rho_0+\sigma) \sigma +
\frac{4}{5}\rho_0^{-\frac{9}{5}}\nabla\rho_0\sigma - h(\sigma,
\rho_0)]\cdot u dx\\
\end{split}
\end{equation}\

Next we estimate the potential part. By the Cauchy-Schwartz
inequality and (5.20), we have
\begin{equation}
\begin{split}
\int(\rho_0+\sigma)^{1+l}\nabla\phi\cdot u dx&\leq \frac{1}{2}
\int(\rho_0+\sigma)^{1+l-\frac{3}{10}}|\nabla\phi|^2dx
+\frac{1}{2} \int(\rho_0+\sigma)^{1+l+\frac{3}{10}}|u|^2dx\\
&\leq C \mathcal{E}_{l+\frac{3}{10}}^0
\end{split}
\end{equation}

The last term in (6.3) is rewritten as
\begin{equation}
\frac{4\pi}{15}\int(\rho_0+\sigma)^{1+l} [
l\text{ }\rho_0^{-\frac{9}{5}}\nabla\rho_0\sigma + \widetilde{h}(\sigma,
\rho_0)]\cdot u dx
\end{equation}
where $\widetilde{h}$ is higher order term including $\sigma$ and
$\nabla\sigma$.  Recall
$\rho_0^{-\frac{9}{5}}\rho_0'=-5r(1+r^2)\sim(\rho_0+\sigma)^
{-\frac{3}{5}}$. When $l=0$, we only have cubic terms left in the
above and hence we are done. The quadratic term when $l\neq 0$ can
be treated as the following:
\begin{equation}
\begin{split}
\int&(\rho_0+\sigma)^{1+l}
[\rho_0^{-\frac{9}{5}}\nabla\rho_0\sigma]\cdot u dx \\
&\leq\frac{1}{2}\int(\rho_0+\sigma)^{1+l+\frac{3}{10}}|u|^2 dx
+\frac{1}{2}\int(\rho_0+\sigma)^{1+l-\frac{3}{10}}
(\rho_0^{-\frac{9}{5}}\rho_0')^2\sigma^2 dx \\
&\leq\frac{1}{2}\int(\rho_0+\sigma)^{1+l+\frac{3}{10}}u^2 dx +
C\int(\rho_0+\sigma)^{-\frac{4}{5}+l+\frac{3}{10}}\sigma^2
dx\\
&\leq C\mathcal{E}_{l+\frac{3}{10}}^0\\
\end{split}
\end{equation}

Apply Lemma 5.1 to the first three integrals in the RHS of
(6.3). With (6.4) and (6.6) the wanted result follows.$\square$\\

For higher order terms, the spirit of details is same as before
but we have extra terms to deal with. While doing higher
derivatives, the necessity of the cooperation with mixed and
spatial estimates occurs.

Let us compute
$$ \int S_l
\partial_t^j \binom{5.2}{5.3} \cdot
\partial_t^j\binom{\sigma}{u} dx=0:$$
\begin{equation}
\begin{split}
\frac{1}{2}&\frac{d}{dt}\{\frac{4\pi}{15}\int(\rho_0+\sigma)^
{-\frac{4}{5}+l}(\partial_t^j\sigma)^2 dx+
\int(\rho_0+\sigma)^{1+l}|\partial_t^j u|^2 dx\} \\
=&\frac{1}{2}\frac{4\pi}{15}(-\frac{4}{5}+l)\int
(\rho_0+\sigma)^{-\frac{4}{5}+l}\frac{\sigma_t}{\rho_0+\sigma}
(\partial_t^j\sigma)^2
dx\\&+\frac{l+1}{2}\int(\rho_0+\sigma)^{1+l}\frac{\sigma_t}
{\rho_0+\sigma}|\partial_t^j u|^2 dx \\
&-\frac{4\pi}{15}\int(\rho_0+\sigma)^{-\frac{4}{5}+l}
\partial_t^j\nabla\cdot[(\rho_0+\sigma)u]\partial_t^j\sigma dx\\
&-\int(\rho_0+\sigma)^{1+l}\partial_t^j[(u\cdot\nabla)
u]\cdot\partial_t^j u
dx\\
&-\frac{4\pi}{15}\int(\rho_0+\sigma)^{1+l}\{\partial_t^j
[(\rho_0+\sigma)^{-\frac{4}{5}} \nabla\sigma] -
\frac{4}{5}\rho_0^{-\frac{9}{5}}\nabla\rho_0\partial_t^j\sigma +
\partial_t^j h(\sigma,
\rho_0)\}\cdot\partial_t^j u dx \\
&-\int(\rho_0+\sigma)^{1+l}\nabla\partial_t^j\phi\cdot\partial_t^j u dx \\
\end{split}
\end{equation}

For computational convenience, we separate $\partial_t^j$ terms
from lower derivative terms in the RHS of (6.7). Some terms contain
unfavorably $(j+1)^{th}$ derivative terms. The worst terms seem to
come from the third and fifth integrals:
$-\frac{4\pi}{15}\int(\rho_0+\sigma)^{-\frac{4}{5}+l}
\nabla\cdot[(\rho_0+\sigma)\partial_t^j u]\partial_t^j\sigma dx$ and
$-\frac{4\pi}{15}\int(\rho_0+\sigma)^{1+l}
[(\rho_0+\sigma)^{-\frac{4}{5}} \nabla\partial_t^j\sigma-
\frac{4}{5}\rho_0^{-\frac{9}{5}}\nabla\rho_0\partial_t^j\sigma]\cdot\partial_t^j
u dx$. Use the integration by parts to get some nice cancellation:
\begin{equation}
\begin{split}
-\frac{4\pi}{15}\int(\rho_0&+\sigma)^{-\frac{4}{5}+l}
\nabla\cdot[(\rho_0+\sigma)\partial_t^j u]\partial_t^j\sigma dx\\
-\frac{4\pi}{15}&\int(\rho_0+\sigma)^{1+l}
[(\rho_0+\sigma)^{-\frac{4}{5}} \nabla\partial_t^j\sigma-
\frac{4}{5}\rho_0^{-\frac{9}{5}}\nabla\rho_0\partial_t^j\sigma]\cdot\partial_t^j
u dx\\
=\frac{4\pi}{15}\int(-\frac{4}{5}&+l)(\rho_0+\sigma)^{-\frac{4}{5}+l}
\nabla(\rho_0+\sigma)\cdot\partial_t^j u\partial_t^j\sigma dx +
(\rho_0+\sigma)^{\frac{1}{5}+l}\partial_t^j
u\cdot\nabla\partial_t^j\sigma dx\\
-\frac{4\pi}{15}&\int(\rho_0+\sigma)^{1+l}
[(\rho_0+\sigma)^{-\frac{4}{5}} \nabla\partial_t^j\sigma-
\frac{4}{5}\rho_0^{-\frac{9}{5}}\nabla\rho_0\partial_t^j\sigma]\cdot\partial_t^j
u dx\\
=\frac{4\pi}{15}\int(\rho_0+&\sigma)^{1+l}[l\text{ }
\rho_0^{-\frac{9}{5}}\nabla\rho_0\partial_t^j\sigma+
\widetilde{h}(\sigma,\nabla\sigma,\rho,\partial_t^j\sigma)]\cdot\partial_t^j
udx
\end{split}
\end{equation}
where
$\widetilde{h}(\sigma,\nabla\sigma,\rho,\partial_t^j\sigma)=(-\frac{4}{5}+l)\{(\rho_0+\sigma)^
{-\frac{9}{5}}\nabla\sigma+[(\rho_0+\sigma)^
{-\frac{9}{5}}-\rho_0^
{-\frac{9}{5}}]\nabla\rho_0\}\partial_t^j\sigma.$  Notice that we
have the above quadratic term only when $l\neq 0$. Taking into account
(6.8) and grouping by similarity, we rewrite (6.7) as following:

\begin{equation}
\begin{split}
\frac{1}{2}&\frac{d}{dt}\{\frac{4\pi}{15}\int(\rho_0+\sigma)^
{-\frac{4}{5}+l}(\partial_t^j\sigma)^2 dx +
\int(\rho_0+\sigma)^{1+l}|\partial_t^j u|^2 dx\} \\
=&\{\frac{1}{2}\frac{4\pi}{15}(-\frac{4}{5}+l)\int
(\rho_0+\sigma)^{-\frac{4}{5}+l}\frac{\sigma_t}{\rho_0+\sigma}
(\partial_t^j\sigma)^2 dx\\
&+\frac{l+1}{2}\int(\rho_0+\sigma)^{1+l}\frac{\sigma_t}
{\rho_0+\sigma}|\partial_t^j u|^2 dx\}\\
+&\{-\frac{4\pi}{15}\int(\rho_0+\sigma)^{-\frac{4}{5}+l}
[\partial_t^j\sigma\nabla\cdot u+\nabla\partial_t^j \sigma\cdot u]
\partial_t^j\sigma dx\\
&-\int(\rho_0+\sigma)^{1+l} [(\partial_t^j u\cdot\nabla) u+
(u\cdot\nabla)\partial_t^j
u]\cdot\partial_t^j u dx\}\\
+&\{\frac{4\pi}{15}\int(\rho_0+\sigma)^{1+l}[l\text{ }
\rho_0^{-\frac{9}{5}}\nabla\rho_0\partial_t^j\sigma+\widetilde{h}
(\sigma,\nabla\sigma,\rho,\partial_t^j\sigma)+
\partial_t^j h(\sigma,
\rho_0)]\cdot\partial_t^j u dx \} \\
+&\{-\int(\rho_0+\sigma)^{1+l}\nabla\partial_t^j\phi\cdot
\partial_t^j u dx\}\\
+&\{-\frac{4\pi}{15}\int(\rho_0+\sigma)^{-\frac{4}{5}+l}
\sum_{i=1}^{j-1}\nabla\cdot[\partial_t^{j-i}(\rho_0+\sigma)\partial_t^i
u]\partial_t^j\sigma dx\\
&-\int(\rho_0+\sigma)^{1+l}\sum_{i=1}^{j-1}[(\partial_t^{i}
u\cdot \nabla)\partial_t^{j-i} u] \cdot\partial_t^j u dx \\
&-\frac{4\pi}{15}\int(\rho_0+\sigma)^{1+l}\sum_{i=1}^{j-1}
[\partial_t^{j-i}(\rho_0+\sigma)^{-\frac{4}{5}}\partial_t^i\nabla\sigma]
\cdot\partial_t^j u dx \}\\
\equiv&(I)+(II)+(III)+(IV)+(V)
\end{split}
\end{equation}\

Note that $(V)$ is alive only for $j=2,\text{ or }3$. First three
groups have the exactly same structure as (6.3). Each term can be
easily estimated as in Lemma 6.2. For $(I)$, by using Lemma 5.2,
we get immediately
\begin{equation}
(I)\leq C\sqrt{\theta_1}\mathcal{E}_l^j\text{ for all }j.
\end{equation}\

Other $(j+1)^{th}$ derivative terms are in $(II)$. After
integrating by parts, the third and fourth integrals become
\begin{equation}
\begin{split}
-\frac{4\pi}{15}&\int(\rho_0+\sigma)^{-\frac{4}{5}+l}
[\partial_t^j\sigma\nabla\cdot u+\nabla\partial_t^j \sigma\cdot u]
\partial_t^j\sigma dx\\
=&-\frac{2\pi}{15}\int(\rho_0+\sigma)^{-\frac{4}{5}+l}(\nabla\cdot
u)(\partial_t^j\sigma)^2dx\\
&+\frac{2\pi}{15}(-\frac{4}{5}+l)\int(\rho_0+\sigma)^{-\frac{4}{5}+l}
\frac{\nabla(\rho_0+\sigma)}{(\rho_0+\sigma)^{\frac{9}{10}}}
\cdot\frac{u}{(\rho_0+\sigma)^{\frac{1}{10}}}
(\partial_t^j\sigma)^2dx\text{,}\\
\
-\int(&\rho_0+\sigma)^{1+l} [(\partial_t^j u\cdot\nabla) u+
(u\cdot\nabla)\partial_t^j
u]\cdot\partial_t^j u dx\\
=&-\int(\rho_0+\sigma)^{1+l} (\partial_t^j u\cdot\nabla)
u\cdot\partial_t^j u
dx+\frac{1}{2}\int(\rho_0+\sigma)^{1+l}(\nabla\cdot u)
|\partial_t^j u|^2dx\\
&+\frac{1+l}{2}\int(\rho_0+\sigma)^{1+l}
\frac{\nabla(\rho_0+\sigma)}{(\rho_0+\sigma)^{\frac{9}{10}}}
\cdot\frac{u}{(\rho_0+\sigma)^{\frac{1}{10}}}|\partial_t^j
u|^2dx\\
\end{split}
\end{equation}\

Therefore (6.11) with Lemma 5.1 gives rise to:
\begin{equation}
(II)\leq C\sqrt{\theta_1}\mathcal{E}_l^j\text{ for all }j.
\end{equation}\

$(III)$ is similar to (6.5) in the zeroth estimate and so it can
be treated in the same way. If we do the same as in (6.6) and use
Lemma 5.1, we get the following: for each $j$,
\begin{equation}
\begin{split}
(III)\leq& C\sqrt{\theta_1} \mathcal{E}_l^j \text{ }\text{ when }l=0,
\text{ and }\\
(III)\leq& C\mathcal{E}_{l+\frac{3}{10}}^j+C\sqrt{\theta_1}
\mathcal{E}_l^j \text{ }\text{ when }l\neq 0.\\
\end{split}
\end{equation}\

The potential part $(IV)$ and cubic terms $(V)$ are somewhat new,
complex and they'd rather be done separately according to
different $j$'s. Before we split the cases, from the dynamics of
$\nabla\partial_t^j\phi$ :
$\nabla\partial_t^j\phi=-4\pi\partial_t^{j-1} [(\rho_0+\sigma)u]$,
we reduce $(IV)$ to the following:
\begin{equation}
(IV)=4\pi\int(\rho_0+\sigma)^{1+l}\partial_t^j u\cdot
\partial_t^{j-1}[(\rho_0+\sigma)u]dx
\end{equation}

Now let $j=1$. Here is the first order $\partial_t$ estimate.

\begin{lem}$(\mathcal{E}_l^1)$ For any small $\eta>0$,
there exist constants $C\text{, }C_{\eta}>0$ such that
\begin{equation*}
\begin{split}
&\frac{1}{2}\frac{d}{dt}\mathcal{E}_0^1\leq
(C\sqrt{\theta_1}+\eta)
\mathcal{E}_0^1+C_{\eta}\mathcal{E}_2^0 \text{ for  }l=0,\\
&\frac{1}{2}\frac{d}{dt}\mathcal{E}_l^1 \leq
(C\sqrt{\theta_1}+\eta) \mathcal{E}_l^1 +
C\mathcal{E}_{l+\frac{3}{10}}^1+ C_{\eta}\mathcal{E}_{l+2}^0
\text{ for }l<0.
\end{split}
\end{equation*}
\end{lem}

Proof. Since $(V)$ has no effect on $j=1$, it is sufficient to
take care of $(IV)$. The potential part is shown to be even better
in terms of derivatives as we can expect in (6.14). By the
Cauchy-Schwartz inequality, we get
\begin{equation}
\begin{split}
(IV)&=4\pi\int(\rho_0+\sigma)^{2+l}\partial_t
u\cdot u dx\\
&\leq \eta\int(\rho_0+\sigma)^{1+l}|\partial_t u|^2 dx+
C_{\eta}\int(\rho_0+\sigma)^{3+l}|u|^2dx\\
&\leq \eta\mathcal{E}_l^1 +C_{\eta}\mathcal{E}_{l+2}^0, \text{ for
any small } \eta>0.\\
\end{split}
\end{equation}

Thus (6.10), (6.12), (6.13) and (6.15) give the desired
result.$\square$\\

The second order $\partial_t^2$ estimate can be done in the same spirit. We have
extra cubic terms to deal with from $(V)$. Lemma 5.2 plays an
important role. Let $j=2$.\

\begin{lem}$(\mathcal{E}_l^2)$ For any small $\eta>0$, there exist
constants $C, C_{\eta}>0$ such that
\begin{equation*}
\begin{split}
&\frac{1}{2}\frac{d}{dt}\mathcal{E}_0^2 \leq
(C\sqrt{\theta_1}+\eta) \mathcal{E}_0^2 +
C_{\eta}\sum_{i=0}^1\mathcal{E}_0^i \text{ for }l=0,\\
&\frac{1}{2}\frac{d}{dt}\mathcal{E}_l^2 \leq
(C\sqrt{\theta_1}+\eta) \mathcal{E}_l^2 +
C\mathcal{E}_{l+\frac{3}{10}}^2+
C_{\eta}\sum_{i=0}^1\mathcal{E}_l^i \text{ for }l<0.
\end{split}
\end{equation*}
\end{lem}\

Proof. The potential part (6.14) can be computed like (6.15).
By Lemma 5.1,
\begin{equation}
\begin{split}
(IV)&=4\pi\int(\rho_0+\sigma)^{2+l}\partial_t^2
u\cdot[\partial_t u + \frac{\sigma_t}{\rho_0+\sigma}u]dx\\
&\leq \eta\mathcal{E}_l^2+C_{\eta}(
\mathcal{E}_{l+2}^1+\theta_1\mathcal{E}_{l+2}^0), \text{ for any small }\eta>0.\\
\end{split}
\end{equation}

Terms in $(V)$ for $j=2$ are at least cubic including mixed
derivatives. We can take the sup for the lowest, first derivative
term and then we end up with manageable quadratic terms. In order
to see how it works, we illustrate the estimate on the first term
in $(V)$:
\begin{equation}
\begin{split}
&\int(\rho_0+\sigma)^{-\frac{4}{5}+l}
\nabla\cdot[\partial_t\sigma\partial_t
u]\partial_t^2\sigma dx\\
=&\int(\rho_0+\sigma)^{\frac{1}{5}+l}
\frac{\partial_t\sigma}{(\rho_0+\sigma)}(\nabla\cdot\partial_t u)
\partial_t^2\sigma dx
+\int(\rho_0+\sigma)^{-\frac{7}{10}+l} \frac{\partial_t u
}{(\rho_0+\sigma)^{\frac{1}{10}}}\cdot\nabla\partial_t\sigma
\partial_t^2\sigma dx\\
\leq&C\sqrt{\theta_1}\{\int(\rho_0+\sigma)^{\frac{6}{5}+l}
|\nabla\cdot\partial_t u|^2
dx+\int(\rho_0+\sigma)^{-\frac{4}{5}+l} (\partial_t^2\sigma)^2
dx\}\\
&+C\sqrt{\theta_1}\{\int(\rho_0+\sigma)^{-\frac{3}{5}+l}
|\nabla\partial_t \sigma|^2
dx+\int(\rho_0+\sigma)^{-\frac{4}{5}+l} (\partial_t^2\sigma)^2
dx\}\\
\leq &C\sqrt{\theta_1}(\widetilde{\mathcal{E}_l}^{1,1}+\mathcal{E}_l^2)\\
\end{split}
\end{equation}

We have used the Cauchy-Schwartz inequality at the first
inequality. Notice the changes of the exponents in weights. We
can do the same to other terms in $(V)$. Ultimately, applying
Lemma 5.2, we get the following:
\begin{equation}
(V)\leq
C\sqrt{\theta_1}\mathcal{E}_l^2+C\sum_{k=0}^{1}\mathcal{E}_
{l+\frac{3}{5}}^k
\end{equation}

Note that $\mathcal{E}_{l+k}\leq C\mathcal{E}_l$ for $k>0$. Thus
(6.10), (6.12), (6.13), (6.16) and (6.18) give the desired
result.$\square$\\

To finish the proof of Proposition 6.1, only $j=3$ i.e. $\partial_t^3$
case is left. The difficulty is to handle the weighted
new cubic terms in $(V)$ of which each factor is at least second
derivative of $\sigma$ and $u$ and hence we cannot utilize Lemma 5.1
directly. To overcome it, we introduce the weighted
Gagliard-Nirenberg inequality. This job is done in the next
lemma. One can see that $(IV)$ and other terms in $(V)$ for $j=3$
can be treated similarly as in (6.16) and (6.17). Therefore, the
following lemma finally establishes Proposition 6.1.\

\begin{lem}Let $l\leq 0$.
For any small fixed $\eta>0$, there exists $C_{\eta}>0$ such that
\begin{equation*}
\begin{split}
&\int(\rho_0+\sigma)^{-\frac{4}{5}+l}
|\partial_t^2\sigma\partial_t\nabla\cdot u\partial_t^3\sigma|
dx,\text{ } \int(\rho_0+\sigma)^{-\frac{4}{5}+l}
|\partial_t\nabla\sigma\cdot\partial_t^2 u\partial_t^3\sigma| dx\\
&\int(\rho_0+\sigma)^{-\frac{4}{5}+l}
|\partial_t^2\sigma\partial_t\nabla\sigma\cdot\partial_t^3u|
dx,\text{ } \int(\rho_0+\sigma)^{1+l} |\partial_t^2
u\cdot\nabla\partial_t u\partial_t^3u| dx\\
\end{split}
\end{equation*}
are bounded by
$\text{ }\eta\mathcal{E}_l^3+C_{\eta}(\mathcal{E}_l^3)^{\frac{3}{2}}
(\mathcal{E}_{l+(-3-2l)}^2)^{\frac{1}{2}}+C_{\eta}
(\mathcal{E}_{l+\frac{1}{5}}^2)^{\frac{3}{2}}
(\mathcal{E}_{l+(-3-2l)}^2)^{\frac{1}{2}}$.\
In particular, if
$\text{ }l\leq-\frac{3}{2}$, they are bounded by
$\eta\mathcal{E}_l+C_{\eta}(\mathcal{E}_l)^2$.
\end{lem}

Proof. First we split each term in
$\partial_t^3$ term and $\partial^2$ term by the Cauchy-Schwartz
inequality. In order to take care of $L^4$ norm of $\partial^2$
 terms we use the
Gagliard-Nirenberg inequality
$\|f\|_{L^4(\mathbb{R}^3)}\leq\frac{1}{2} \|\nabla
f\|_{L^2(\mathbb{R}^3)}^{\frac{3}{4}} \|f\|_{L^2(\mathbb{R}^3)}
^{\frac{1}{4}}$. Since the inequality complies well with the
localization and the weight functions are nice, using a partition
of unity, in our case we get the weighed version of the
Gagliard-Nirenberg inequality

\begin{equation}
\begin{split}
\int_0^{\infty}w_k f^4 r^2 dr \leq C&
(\int_0^{\infty}w_k^{\alpha}|\nabla f|^2 r^2 dr)^{\frac{3}{2}}
(\int_0^{\infty}w_k^{\beta} f^2 r^2 dr) ^{\frac{1}{2}}\\ +&C
(\int_0^{\infty}w_k^{\alpha'}f^2 r^2 dr)^{\frac{3}{2}}
(\int_0^{\infty}w_k^{\beta'} f^2 r^2 dr) ^{\frac{1}{2}}
\end{split}
\end{equation}
where $w_k=(1+r^2)^{\frac{k-5}{2}}\sim
(\rho_0+\sigma)^{1-\frac{k}{5}} $ and
$\frac{3}{2}\alpha+\frac{1}{2}\beta=
\frac{3}{2}\alpha'+\frac{1}{2}\beta'=1$. Its proof
is given at the end of the argument. Only the first and the
last term are treated in this proof. The other cases can be
estimated in the same way.

\begin{equation}
\begin{split}
\int_0^{\infty}(\rho_0&+\sigma)^{-\frac{4}{5}+l}
|\partial_t^2\sigma\partial_t\nabla\cdot u\partial_t^3\sigma|r^2 dr \\
\leq &\eta
\int_0^{\infty}(\rho_0+\sigma)^{-\frac{4}{5}+l}(\partial_t^3\sigma)^2
r^2 dr +C_{\eta}\int_0^{\infty}(\rho_0+\sigma)^{-\frac{4}{5}+l}
|\partial_t^2\sigma\partial_t \nabla\cdot u|^2 r^2 dr\\
\
\int_0^{\infty}(\rho_0&+\sigma)^{-\frac{4}{5}+l}
|\partial_t^2\sigma\partial_t \nabla\cdot u|^2 r^2 dr \\
\leq&\frac{1}{2}\int_0^{\infty}(\rho_0+\sigma)^{-\frac{14}{5}+l}
(\partial_t^2\sigma)^4 r^2
dr+\frac{1}{2}\int_0^{\infty}(\rho_0+\sigma)^{\frac{6}{5}+l}
(\partial_t \nabla\cdot u)^4 r^2 dr\\
\end{split}
\end{equation}

\begin{equation}
\begin{split}
\int_0^{\infty}(\rho_0&+\sigma)^{-\frac{14}{5}+l}
(\partial_t^2\sigma)^4 r^2 dr\\
 \leq &C \{\int_0^{\infty}(\rho_0+\sigma)^{-\frac{3}{5}+l}
|\partial_t^2\nabla\sigma|^2 r^2
dr\}^{\frac{3}{2}}\{\int_0^{\infty}
(\rho_0+\sigma)^{-\frac{19}{5}-l} (\partial_t^2\sigma)^2 r^2
dr\}^{\frac{1}{2}}\\
&+C\{\int_0^{\infty}(\rho_0+\sigma)^{-\frac{3}{5}+l}
(\partial_t^2\sigma)^2 r^2 dr\}^{\frac{3}{2}}\{\int_0^{\infty}
(\rho_0+\sigma)^{-\frac{19}{5}-l} (\partial_t^2\sigma)^2 r^2
dr\}^{\frac{1}{2}}\\
 \leq& C(\mathcal{E}_l^3)^{\frac{3}{2}}(\mathcal{E}_{l+(-3-2l)}^2)
^{\frac{1}{2}}+C(\mathcal{E}_{l+\frac{1}{5}}^2)^{\frac{3}{2}}
(\mathcal{E}_{l+(-3-2l)}^2)^{\frac{1}{2}}\\
\
\int_0^{\infty}(\rho_0&+\sigma)^{\frac{6}{5}+l}|\partial_t\nabla
\cdot u|^4
r^2 dr\\
 \leq &C \{\int_0^{\infty}(\rho_0+\sigma)^{\frac{7}{5}+l}
|\partial_t\nabla(\nabla\cdot u)|^2 r^2
dr\}^{\frac{3}{2}}\{\int_0^{\infty}
(\rho_0+\sigma)^{-\frac{9}{5}-l} |\partial_t \nabla\cdot u|^2 r^2
dr\}^{\frac{1}{2}}\\
&+C \{\int_0^{\infty}(\rho_0+\sigma)^{\frac{7}{5}+l}
|\partial_t\nabla\cdot u|^2 r^2
dr\}^{\frac{3}{2}}\{\int_0^{\infty}
(\rho_0+\sigma)^{-\frac{9}{5}-l} |\partial_t \nabla\cdot u|^2 r^2
dr\}^{\frac{1}{2}}\\
\leq& C(\mathcal{E}_l^3)^{\frac{3}{2}}(\mathcal{E}_{l+(-3-2l)}^2)^
{\frac{1}{2}}+C(\mathcal{E}_{l+\frac{1}{5}}^2)^{\frac{3}{2}}
(\mathcal{E}_{l+(-3-2l)}^2)^{\frac{1}{2}}\\
\end{split}
\end{equation}\

Here is the verification of each exponent:
$-\frac{14}{5}+l=\frac{3}{2}(-\frac{3}{5}+l)
+\frac{1}{2}(-\frac{19}{5}-l)$ and
$-\frac{19}{5}-l=-\frac{4}{5}+l+(-3-2l) $;
$\frac{6}{5}+l=\frac{3}{2}(\frac{7}{5}+l)+\frac{1}{2}
(-\frac{9}{5}-l)$ and $-\frac{9}{5}-l=\frac{6}{5}+l+(-3-2l)$. Note
that $-\frac{19}{5}-l\geq-\frac{4}{5}+l$ for $l\leq-\frac{3}{2}$.
Now let us look at the last term. We go through the similar
computation as in (6.20) and (6.21).

\begin{equation}
\begin{split}
\int_0^{\infty}(\rho_0+ \sigma&)^{1+l}|\partial_t^2
u\cdot\nabla\partial_t u\partial_t^3u| r^2 dr \leq \eta
\int_0^{\infty}(\rho_0+ \sigma)^{1+l}|\partial_t^3 u|^2 r^2 dr \\
&+C_{\eta}\{\int_0^{\infty}(\rho_0+\sigma)^{\frac{4}{5}+l}
|\partial_t^2 u|^4 r^2 dr +\int_0^{\infty}(\rho_0+
\sigma)^{\frac{6}{5}+l}
|\nabla\partial_t u|^4 r^2 dr\} \\
\end{split}
\end{equation}

\begin{equation}
\begin{split}
\int_0^{\infty}(\rho_0&+\sigma)^{\frac{4}{5}+l}|\partial_t^2 u|^4
r^2 dr \\ \leq &C \{\int_0^{\infty}(\rho_0+\sigma)^{\frac{6}{5}+l}
|\partial_t^2 \nabla u|^2 r^2 dr\}^{\frac{3}{2}}\{\int_0^{\infty}
(\rho_0+\sigma)^{-2-l} |\partial_t^2 u|^2 r^2 dr\}^{\frac{1}{2}}\\
&+C \{\int_0^{\infty}(\rho_0+\sigma)^{\frac{6}{5}+l}
|\partial_t^2 u|^2 r^2 dr\}^{\frac{3}{2}}\{\int_0^{\infty}
(\rho_0+\sigma)^{-2-l} |\partial_t^2 u|^2 r^2 dr\}^{\frac{1}{2}}\\
 \leq &C
(\mathcal{E}_l^3)^{\frac{3}{2}}(\mathcal{E}_{l+(-3-2l)}^2)
^{\frac{1}{2}}+C
(\mathcal{E}_{l+\frac{1}{5}}^2)^{\frac{3}{2}}(\mathcal{E}_{l+(-3-2l)}^2)
^{\frac{1}{2}}\\
\
\int_0^{\infty}(\rho_0&+\sigma)^{\frac{6}{5}+l}|\nabla\partial_t
u|^4 r^2 dr \\ \leq &C
(\mathcal{E}_l^3)^{\frac{3}{2}}(\mathcal{E}_{l+(-3-2l)}^2)
^{\frac{1}{2}}+C
(\mathcal{E}_{l+\frac{1}{5}}^2)^{\frac{3}{2}}(\mathcal{E}_{l+(-3-2l)}^2)
^{\frac{1}{2}}, \text{ from (6.21)}\\
\end{split}
\end{equation}\

Observe that
$\frac{4}{5}+l=\frac{3}{2}(\frac{6}{5}+l)+\frac{1}{2}(-2-l)$,
$\frac{6}{5}+l=\frac{3}{2}(\frac{7}{5}+l)+\frac{1}{2}(-\frac{9}{5}-l)$
and $-2-l=1+l+(-3-2l)$, $-\frac{9}{5}-l=\frac{6}{5}+l+(-3-2l)$.

The only missing part is the proof of the weighted
Gagliard-Nirenberg inequality. Here it comes:

Proof of (6.19): We choose a partition of unity
$\{\varphi_n\}_{n\geq 0}$ as following:
\[
\varphi_0(r)=
  \begin{cases}
    1-r & \text{if $0\leq r \leq 1$} \\
    0 & \text{if $r\geq\ 1$}\\
  \end{cases}
\]
\[
\text{ }\text{ }\text{ }\text{ }\text{ }\text{ }\text{ }\text{
}\text{ }\text{ }\text{ }\varphi_n(r)=
  \begin{cases}
    0 & \text{if $0\leq r \leq \frac{n-1}{2}$} \\
    r-\frac{n-1}{2} & \text{if $\frac{n-1}{2}\leq r \leq
    \frac{n}{2}$} \\
    \frac{1}{2} & \text{if $\frac{n}{2}\leq r \leq \frac{n+1}{2}$}\\
    \frac{n+2}{2}-r & \text{if $\frac{n+1}{2}\leq r\leq \frac{n+2}{2}$}\\
    0 &\text{if $r\geq \frac{n+2}{2}$}\\
  \end{cases}
\]

It is easy to check $0\leq \varphi_n\leq 1$, $\text{supp
}\varphi_n=[\frac{n-1}{2},\frac{n+2}{2}]$( $\text{supp
}\varphi_0=[0,1]$), $\sum_{n=0}^{\infty}\varphi_n(r)=1\text{ for
all }r\geq 0$, and $|\varphi_n'(r)|\leq 1$ a.e.
\begin{equation}
\begin{split}
\int_0^{\infty}(1+r^2)^k
f^4r^2dr=&\int_0^{\infty}[\sum_{n=0}^{\infty}\varphi_n(r)]
^4(1+r^2)^k
f^4r^2dr\\
\leq &C\sum_{n=0}^{\infty}\int_0^{\infty}\varphi_n^4(r)(1+r^2)^k
f^4r^2dr\\
\end{split}
\end{equation}

We only consider $k\geq 0$. Other cases can be proven in the same
manner. Note that $(1+r^2)^k\leq(1+(\frac{n+2}{2})^2)^k$ on
$\text{supp }\varphi_n$. First we localize the half real line
according to the partition of unity. Since weights are monotonic,
they can be localized as well. And then we apply the
Gagliardo-Nirenberg inequality.

\begin{equation*}
\begin{split}
(9.1)\leq&C\sum(1+(\frac{n+2}{2})^2)^k\int_0^{\infty}\varphi_n^4
f^4 r^2 dr\\ \leq&
C\sum(1+(\frac{n+2}{2})^2)^k(\int_0^{\infty}[(\varphi_n f)' ]^2
r^2 dr)^{\frac{3}{2}} (\int_0^{\infty}{\varphi_n}^{2}f^2
r^2 dr )^{\frac{1}{2}}\\
\leq& C\sum(1+(\frac{n+2}{2})^2)^k(\int_0^{\infty} [(\varphi_n f'
)^2+(\varphi_n' f )^2]r^2 dr)^{\frac{3}{2}}
(\int_0^{\infty}{\varphi_n}^{2}f^2 r^2 dr )^{\frac{1}{2}}\\
\leq&
C\sum(1+(\frac{n+2}{2})^2)^k(\int_{\frac{n-1}{2}}^{\frac{n+2}{2}}
 f'^2r^2 dr)^{\frac{3}{2}}
(\int_{\frac{n-1}{2}}^{\frac{n+2}{2}}f^2 r^2 dr )^{\frac{1}{2}}\\
&+C\sum(1+(\frac{n+2}{2})^2)^k(\{\int_{\frac{n-1}{2}}^{\frac{n}{2}}
+\int_{\frac{n+1}{2}}^{\frac{n+2}{2}}\} f^2r^2 dr)^{\frac{3}{2}}
(\int_{\frac{n-1}{2}}^{\frac{n+2}{2}}f^2 r^2 dr )^{\frac{1}{2}}\\
\leq&C\sum(\int_{\frac{n-1}{2}}^{\frac{n+2}{2}}
 (1+r^2)^{\alpha}f'^2r^2 dr)^{\frac{3}{2}}
(\int_{\frac{n-1}{2}}^{\frac{n+2}{2}}(1+r^2)^{\beta}f^2 r^2 dr
)^{\frac{1}{2}}\\
&+C\sum(\{\int_{\frac{n-1}{2}}^{\frac{n}{2}}
+\int_{\frac{n+1}{2}}^{\frac{n+2}{2}}\}(1+r^2)^{\alpha'} f^2r^2
dr)^{\frac{3}{2}}
(\int_{\frac{n-1}{2}}^{\frac{n+2}{2}}(1+r^2)^{\beta'}f^2 r^2 dr )
^{\frac{1}{2}}\\
\leq&C(\int_0^{\infty}(1+r^2)^{\alpha}f'^2r^2 dr
)^{\frac{3}{2}}(\int_0^{\infty}(1+r^2)^{\beta}f^2 r^2 dr
)^{\frac{1}{2}}\\ &+C(\int_0^{\infty}(1+r^2)^{\alpha'}f^2r^2 dr
)^{\frac{3}{2}}(\int_0^{\infty}(1+r^2)^{\beta'}f^2 r^2 dr
)^{\frac{1}{2}}\\
\end{split}
\end{equation*}
In the above $C$ is a generic constant. Note that
$\frac{3}{2}\alpha+\frac{1}{2}\beta=\frac{3}{2}\alpha'
+\frac{1}{2}\beta'=k$.$\square$\\

\section{Nonlinear Instability}
Now we are ready to prove the bootstrap argument. The proof
completely depends on the estimates in Section 6 and the Gronwall
inequality.

\begin{prop} Let $\nu(t)=\binom{\sigma(t)}{u(t)}$ be a solution of
the Euler-Poisson system (1.9) and (1.10). Let $l_{\star}\leq -3$
be given. Assume that
$$\sqrt{\mathcal{E}_{l_{\star}}}(0)\leq
C_0\delta \text{ and } \sqrt{\mathcal{E}_0^0}(t)\leq C_0 \delta
e^{\sqrt{\mu_0} t} \text{ for } 0\leq t \leq T.$$ Then there
exist $C_5$, $\theta_0$ $>0$ such that
$$\text{if } \text{ }
0\leq t \leq \min\{T, T^{\delta}\}, \text{ then }
\sqrt{\widetilde{\mathcal{E}_{l_{\star}}}}(t) \leq C_5\delta
e^{\sqrt{\mu_0} t} \leq C_5 \theta_0,$$ where $T^{\delta}=
\frac{1}{\sqrt{\mu_0}}\ln \frac{\theta_0}{\delta}.$
\end{prop}

Proof. In Proposition 6.1, Lemma 6.2, 6.3 and 6.4, choose
$\theta_1$ and $\eta$ small enough so that
$C\sqrt{\theta_1}+\eta\leq \frac{\sqrt{\mu_0}}{2}$. Therefore
there exist constants $C_1\geq 0$, $C_2, C_3>0$ such that for
$l\leq 0$,

\begin{equation*}
\begin{split}
&(a_l)\text{ } \frac{1}{2}\frac{d}{dt}\mathcal{E}_l^0 \leq
\frac{\sqrt{\mu_0}}{2}
\mathcal{E}_l^0+C_2\mathcal{E}_{l+\frac{3}{10}}^0,\\
&(b_l)\text{ } \frac{1}{2}\frac{d}{dt}\mathcal{E}_l^1 \leq
\frac{\sqrt{\mu_0}}{2}
\mathcal{E}_l^1+C_1\mathcal{E}_{l+\frac{3}{10}}^1 +C_2\mathcal{E}
_{l+2}^0,\\
&(c_l)\text{ } \frac{1}{2}\frac{d}{dt}\mathcal{E}_l^2 \leq
\frac{\sqrt{\mu_0}}{2}\mathcal{E}_l^2+C_1\mathcal{E}_{l+\frac{3}{10}}^2
+
C_2(\mathcal{E}_l^0+\mathcal{E}_l^1),\\
&(d_l)\text{ }
\frac{1}{2}\frac{d}{dt}\mathcal{E}_l\leq\frac{\sqrt{\mu_0}}{2}\mathcal{E}_l+
C_1\mathcal{E}_{l+\frac{3}{10}}
+C_2(\mathcal{E}_l^0+\mathcal{E}_l^1+\mathcal{E}_l^2)+
C_3(\mathcal{E}_l)^{\frac{3}{2}}(\mathcal{E}_{l+(-3-2l)})^{\frac{1}{2}}.
\end{split}
\end{equation*}\

Note that $C_1=0$ when $l=0$. Define $T^{\ast}$ by
$$T^{\ast}\equiv \sup\{t:\widetilde{\mathcal{E}_{l}}(s)\leq \min\{\theta_1,
\frac{\sqrt{\mu_0}}{4C_3}\}, s\in[0,t]\text{, }l_\star \leq l\leq
0\}.$$

Here $\theta_1$ is a small constant coming from the smallness
assumption to guarantee that the nonlinear energy estimates work.

Let $0\leq t \leq \min\{T,T^{\ast}\}$. Then since
$\sqrt{\mathcal{E}_0^0}(t)\leq C_0 \delta e^{\sqrt{\mu_0} t}$
by the hypothesis,
from $(a_l)$, we get
$$(a_l)\Longrightarrow\sqrt{\mathcal{E}_l^0}(t)\leq C_0' \delta
e^{\sqrt{\mu_0} t} \text{ for all }l\leq 0$$ by the standard
Gronwall inequality. We use $C_0'$ as a generic constant.
Consider the following diagram:
\begin{equation*}
\begin{split}
(b_0)\text{ }\text{ }&\text{ }\frac{1}{2}
\frac{d}{dt}\mathcal{E}_0^1 \leq
\frac{\sqrt{\mu_0}}{2}
\mathcal{E}_0^1+C_2\mathcal{E}_2^0\Longrightarrow
\sqrt{\mathcal{E}_0^1}(t)\leq C_0' \delta e^{\sqrt{\mu_0} t}\\
(b_{-\frac{3}{10}})&\text{ }\frac{1}{2}\frac{d}{dt}\mathcal{E}_{-\frac{3}{10}}^1 \leq
\frac{\sqrt{\mu_0}}{2}
\mathcal{E}_{-\frac{3}{10}}^1+C_1\mathcal{E}_0^1 +C_2\mathcal{E}
_{\frac{17}{10}}^0\Longrightarrow\sqrt{\mathcal{E}_{-\frac{3}{10}}^1}(t)\leq
C_0' \delta e^{\sqrt{\mu_0} t}\\
\end{split}
\end{equation*}

Likewise, for all $k\geq 0$, we have
\begin{equation*}
(b_{{-\frac{3}{10}}k})\Longrightarrow
\sqrt{\mathcal{E}_{-\frac{3}{10}k}^1}(t)\leq C_0' \delta
e^{\sqrt{\mu_0} t}.
\end{equation*}

By the same bootstrap argument with $(c_{{-\frac{3}{10}}k})$
where $k\geq 0$, we have
$$\sqrt{\mathcal{E}_{-\frac{3}{10}k}^2}(t)\leq C_0' \delta
e^{\sqrt{\mu_0} t}.$$

Now we move onto $(d_l)$. Let us start with $l=0$. Recall that
$C_1=0$. Since $(\mathcal{E}_{0})^{\frac{1}{2}}
(\mathcal{E}_{-3})^{\frac{1}{2}}\leq\frac{\sqrt{\mu_0}}{4C_3}$, firstly we
get

\begin{equation*}
\begin{split}
\frac{d}{dt}\mathcal{E}_0&\leq\sqrt{\mu_0}\mathcal{E}_0+2C_2
(\mathcal{E}_0^0+\mathcal{E}_0^1+\mathcal{E}_0^2)+
2C_3(\mathcal{E}_0)^{\frac{3}{2}}(\mathcal{E}_{-3})^{\frac{1}{2}}\\
&\leq\frac{3\sqrt{\mu_0}}{2} \mathcal{E}_0 + 2C_2C_0'^2
\delta^2e^{2\sqrt{\mu_0}t}
\end{split}
\end{equation*}\

By the  Gronwall inequality, we have the exponential growth with
the exponent $2\sqrt{\mu_0}$ on $\mathcal{E}_0$. As the previous
cases, form a sequence $\{(d_{-\frac{3}{10}k})\}_{k\geq 0}$ and
use the Gronwall inequality to get the following
\begin{equation}
\mathcal{E}_l(t) \leq C_4^2\delta^2e^{2\sqrt{\mu_0}t},\text{ for }
0\leq t \leq \min\{T,T^{\ast}\}
\end{equation}
for any $l\leq 0$ and some constant $C_4$. And in success by
Lemma 5.2, for any $l\leq 0$, we also have
\begin{equation}
\widetilde{\mathcal{E}_l}(t) \leq
C_5^2\delta^2e^{2\sqrt{\mu_0}t},\text{ for } 0\leq t \leq
\min\{T,T^{\ast}\}.
\end{equation}\

Now choose $\theta_0$ such that $(C_5\theta_0)^2 <
\min\{\theta_1, \frac{\sqrt{\mu_0}}{4C_3}\}$. We consider
following two cases:\\

(i) $T^{\delta}\leq \min\{T,T^{\ast}\}$; in this case, the
conclusion follows without any extra work.\\

(ii) $T^{\delta}> \min\{T,T^{\ast}\}$; then $T\leq T^{\ast}
<T^{\delta}$. If this is true, then again the conclusion is
trivial. We show this is the only possibility. If not, we have\\
$T^{\ast}< T <T^{\delta}$. Letting $t=T^{\ast}$, from (7.2) and
the definition of $T^{\delta}$, we get
$$\widetilde{\mathcal{E}_l}(T^{\ast})\leq C_5^2\delta^2e^{2\sqrt{\mu_0}T^{\ast}}<
(C_5\delta \exp^{\sqrt{\mu_0}T^{\delta}})^2=(C_5 \theta_0)^2.$$
But this is impossible by the choice of $\theta_0$ since it
would contradict the definition of $T^{\ast}$. This establishes
the proposition.$\square$\\

From now on we regard $\delta>0$ as an arbitrary small parameter
and $\theta$ as a small but fixed positive constant independent
of $\delta$. Recall that $T^{\delta}$ is defined by $\theta=\delta
\exp^{\sqrt{\mu_0}T^{\delta}}$ or equivalently $T^{\delta}=
\frac{1}{\sqrt{\mu_0}}\ln
\frac{\theta}{\delta}$.\\

$\mathbf{Proof\text{ }of\text{ }Theorem\text{ }1.1.}$ Let $\nu_0
= \binom{\phi_0}{\psi_0}$ be a growing mode for the linearized
Euler-Poisson system (2.1) and (2.2) obtained in Section 2.
Normalize $\nu_0$ such that $$\|\nu_0\|_0^2\equiv\frac{16}{15}
\pi^2(\|\phi_0\|_{V_0}^2
+\|\psi_0\|_{W_0}^2)=1.$$ By the behavior
of $\psi_0$ in Proposition 3.4 and the relation (2.4) between $\phi_0$
and $\psi_0$, we know $|\delta\nu_0|$ is relatively
small compared to $\rho_0$; $\rho_0 +\delta\phi_0 \sim \gamma\rho_0$
where $\gamma$ is close to 1, and $\int_0^{\infty}\phi_0 r^2 dr=0$. We may assume $\frac{19}{20}\rho_0
< \rho_0 +\delta\phi_0 < \frac{21}{20}\rho_0$. Let
$$4\pi\int_0^{\infty}\frac{4\pi}{15}
(\rho_0+\delta\phi_0)^{-\frac{4}{5}+l}\phi_0^2 r^2 dr +4\pi
\int_0^{\infty}(\rho_0+\delta\phi_0)^{1+l}\psi_0^2 r^2 dr = a^2 <
\infty.$$ Now solve the Euler-Poisson system with a family of
initial data $\nu|_{t=0}=\delta \nu_0$. The continuity equation
gives rise to $\int_0^{\infty}\sigma^{\delta}r^2dr=0$. Denote the corresponding
$\mathcal{E}_l$-solution by
$\nu(t)\equiv\nu^{\delta}(t)=\binom{\sigma^{\delta}(t)}{u^{\delta}(t)}$.
It can be written as
$$\nu(t)=\delta e^{\sqrt{\mu_0}t}\nu_0 +
\int_0^t\mathcal{L}(t-\tau)N(\tau)d\tau.$$ where $\mathcal{L}$ is
the solution operator for the linearized Euler-Poisson system and
$N$ is nonlinear part.
$$N=\binom{\frac{1}{r^2}(r^2\sigma u)_r}{uu_r+\frac{4\pi}{15}
[-\frac{4}{5}\rho_0^{-\frac{9}{5}}\sigma\sigma_r
+(\rho_0+\sigma)_r h]}$$ where $h=(\rho_0+\sigma)^{-\frac{4}{5}}
-\rho_0^{-\frac{4}{5}}+\frac{4}{5}\rho_0^{-\frac{9}{5}}\sigma$.\

Define $T$ by
$$T\equiv\sup\{s:\text{ for }0\leq t \leq s, \text{ }
\sqrt{\mathcal{E}_0^0}\leq \max\{3, a\}\delta
e^{\sqrt{\mu_0}t}\}.$$

Then by Proposition 7.1, there exist $C_1$ and $\theta_0$ $>0$
such that \\for $0\leq t \leq \min\{T, T^{\delta}\}$,
$$\sqrt{\mathcal{E}_l}(t)\leq\sqrt{\widetilde{\mathcal{E}_l}}(t)
\leq C_1 \delta e^{\sqrt{\mu_0}t} \leq
C_1\theta_0.$$

Note for sufficiently small $\theta_0$, by Lemma 5.1, it means
that $|\frac{\sigma^{\delta}}{\rho_0+\sigma^{\delta}}|<< 1$ i.e. $
\rho_0+\sigma^{\delta}$ behaves like $\rho_0$. So we can find a
small constant $\beta$, $0 \leq \beta \leq \frac{1}{2}$ such that
\begin{equation}
\begin{split}
(1-\beta)^2\rho_0^{-\frac{4}{5}}& \leq
(\rho_0+\sigma^{\delta})^{-\frac{4}{5}} \leq
(1+\beta)^2\rho_0^{-\frac{4}{5}}\\
(1-\beta)^2\rho_0 &\leq \rho_0+\sigma^{\delta} \leq
(1+\beta)^2\rho_0.
\end{split}
\end{equation}

To emphasize which functions we deal with, we denote
$\mathcal{E}_0^0$ being plugged $f=\binom{f_1}{f_2}$ in
but its weight part unchanged by
$\|f\|_Y^2$:
$$\|f\|_Y^2=4\pi\int_0^{\infty}\frac{4\pi}{15}(\rho_0+\sigma^{\delta})
^{-\frac{4}{5}}(f_1)^2r^2dr+4\pi\int_0^{\infty}(\rho_0+\sigma^{\delta})
(f_2)^2r^2dr $$

In this notation, by (7.3), for $0\leq t \leq
\min\{T, T^{\delta}\}$,
\begin{equation}
\|\delta e^{\sqrt{\mu_0}t}\nu_0 \|_Y \geq (1-\beta)\delta
e^{\sqrt{\mu_0}t}\|\nu_0\|_0= (1-\beta)\delta
e^{\sqrt{\mu_0}t}.
\end{equation}\

For the nonlinear parts, from the linearized estimates in the
Lemma 4.1, 4.2 and 4.3, for $0\leq t \leq \min\{T, T^{\delta}\}$,
we have
\begin{equation}
\begin{split}
\|\nu^{\delta}(t)-\delta e^{\sqrt{\mu_0}t}\nu_0
\|_{Y}=&\|\int_0^t\mathcal{L}(t-\tau)N(\tau)d\tau
\|_{Y}\\
\leq &C\int_0^t
e^{\sqrt{\mu_0}(t-\tau)}(\|N(\tau)\|_Y+\|\partial_t N(\tau)\|_Y
+\|\partial_t^2 N(\tau)\|_Y) d\tau\\
\leq &C\int_0^te^{\sqrt{\mu_0}(t-\tau)} \{\text{
}|\frac{\sigma}{\rho_0+\sigma}|+ |\frac{\sigma_t}{\rho_0+\sigma}|+
|\frac{\nabla\sigma}{(\rho_0+\sigma) ^{\frac{9}{10}}}|\\
&\text{ }+|\frac{u}{(\rho_0+\sigma)^{\frac{1}{10}}}|+
|\frac{u_t}{(\rho_0+\sigma)^{\frac{1}{10}}}|+|\nabla u|\text{
}\}\text{ }
\{\sqrt{\mathcal{E}}_l+\mathcal{E}_l\}d\tau\\
\leq& C\int_0^te^{\sqrt{\mu_0}(t-\tau)}(\delta
e^{\sqrt{\mu_0}\tau})
(\delta e^{\sqrt{\mu_0}\tau})d\tau\\
\leq& C_2(\delta e^{\sqrt{\mu_0}t})^2
\end{split}
\end{equation}
where $C_2$ is a constant. At the second inequality we have used
Lemma 5.2. The next inequality follows from Proposition 7.1 and
Lemma 5.1.\

Now if necessary, fix $\theta_0$ sufficiently small such that
$C_2\theta_0 \leq \frac{1-\beta}{2}$.\\

Claim. $T^{\delta}\leq T$.\\

Proof. If not i.e. $T^{\delta} > T$, by (7.4) and (7.5)

\begin{equation*}
\begin{split}
\|\nu^{\delta}\|_{Y}(T) &\leq \|\delta e^{\sqrt{\mu_0}t}\nu_0
\|_{Y}(T)+ \|\nu^{\delta}-\delta e^{\sqrt{\mu_0}t}\nu_0
\|_{Y}(T)\\
&\leq(1+\beta)\delta e^{\sqrt{\mu_0}T}\|\nu_0\|_0 +
C_2\theta_0
\delta e^{\sqrt{\mu_0}T}\\
&\leq \frac{3+\beta}{2}\delta e^{\sqrt{\mu_0}T} < 2\delta
e^{\sqrt{\mu_0}T}\\
\end{split}
\end{equation*}
which would contradict the definition of $T$.\\

Once we have $T^{\delta}\leq T$, again by (7.4) and (7.5),

$$\sqrt{\mathcal{E}_0^0}(T^{\delta}) \geq (1-\beta)
\delta e^{\sqrt{\mu_0}T^{\delta}}-\frac{1-\beta}{2}\delta
e^{\sqrt{\mu_0}T^{\delta}}=\frac{1-\beta}{2}\theta_0 >0.$$

Set $\theta=\frac{1-\beta}{2}\theta_0$. This finishes the
proof of the theorem.$\square$\\

\textbf{Acknowledgments:} The author would like to deeply
thank \textsc{Yan Guo}
for many stimulating feedbacks and inspiring discussions.\\


\end{document}